\documentclass[11pt]{amsart}

\calclayout

\usepackage{mathrsfs}

\usepackage{amssymb}

\usepackage[normalem]{ulem}

\usepackage{color}

\usepackage[fixamsmath]{mathtools}

\usepackage{tikz-cd}

\newtheorem{proposition}{Proposition}

\theoremstyle{definition}

\theoremstyle{remark}

\renewcommand{\Re}{\hp\mbox{Re}\pt} 
\renewcommand{\Im }{\hp\mbox{Im}\hp}

\newcommand{\pt}{\hspace{1pt}} 
\newcommand{\hp}{\hspace{0.5pt}} 

\newcommand{\I}{{\scriptscriptstyle I}} 
\newcommand{\J}{{\scriptscriptstyle J}}

  {\list{}{\leftmargin=0.8in\rightmargin=0in}\item[]}%
  {\endlist}

\begin{document}


\title[Gravitational instantons of type $D_k$]{Gravitational instantons of type $D_k$ and a generalization of the Gibbons-Hawking Ansatz}


\author{Radu A. Iona\c{s}}

\newcommand{\Addresses}{\noindent \textsc{C. N. Yang Institute for Theoretical Physics} \\ \textsc{Stony Brook University, Stony Brook, NY 11794, U.S.A.}}

%
%

\date{}



\thispagestyle{empty}



\begin{abstract}
We describe a quaternionic-based Ansatz generalizing the Gibbons-Hawking Ansatz to a class of hyperk\"ahler metrics with hidden symmetries. We then apply it to obtain explicit expressions for gravitational instanton metrics of type $D_k$. 
\end{abstract}

\maketitle

\vskip70pt

\section{Introduction}

Through the guise of hypercomplex structures, quaternions play a central role in the definition of hyperk\"ahler spaces. However, in practice, as one progresses towards more explicit descriptions of hyperk\"ahler metrics, the quaternionic structure is, as a rule, relegated to relative obscurity while some symmetry principle takes the center stage instead. Thus is the case for example of the class of hyperk\"ahler metrics described by the Gibbons-Hawking Ansatz, where  an abelian tri-Hamiltonian symmetry plays the leading role. Or that of the Atiyah-Hitchin metric on the centered moduli space of charge-2 monopoles, to give another well-known example, where the manifest role is played by an $SO(3)$ symmetry rotating the hyperk\"ahler structure. 

In this paper we will try to make the case that this is not as much a necessity as a choice. We will argue that it is, in other words, possible to bring the quaternions up front\,---\,in the form of quaternionic coordinates for the hyperk\"ahler space\,---\,and relegate the symmetries to the background. While for classical symmetries this offers nothing more than an alternative, the real advantage of this approach becomes apparent when we consider hyperk\"ahler metrics with \textit{hidden} symmetries, for which we do not have entirely satisfactory symmetry-based descriptions. 

More specifically, we do two things. First, we recast a class of hyperk\"ahler metrics with hidden symmetries discovered by Lindstr\"om and Ro\v{c}ek into a quaternionic coordinate frame to obtain an Ansatz which for this class is essentially what the Gibbons-Hawking Ansatz is for the class of hyperk\"ahler metrics with maximal toric tri-Hamiltonian symmetry. The key to achieving this consists in expressing a certain spinorial object in spherical (read \textit{quaternionic}) coordinates, a process which naturally requires the introduction of an auxiliary elliptic curve. Second, we use this Ansatz to give an explicit description of gravitational instanton metrics of type $D_k$ in terms of associated quasi-elliptic functions. 

Gravitational instantons are by definition complete connected four-dimensional hyper\-k\"ahler Riemannian manifolds. When the manifold is non-compact with one end this definition is supplemented with a requirement that the metric approaches according to a certain prescription the flat metric at infinity. Thus, \textit{asymptotically locally Euclidean} (ALE) metrics are required to approach the flat metric on $\mathbb{R}^4/\Gamma$, with $\Gamma$ a finite subgroup of $SU(2)$, at a certain decay rate (as well as appropriate decay rates in the derivatives).\footnote{In this case, incidentally, the manifold can only have one end, and so the notion of its ``infinity" is unambiguous in this respect \cite{MR3543182}.} According to whether $\Gamma$ is a cyclic, binary dihedral, tetrahedral, octahedral or dodecahedral subgroup of $SU(2)$, these are divided into two infinite discrete series, $A_k$ and $D_k$, and three exceptional types, $E_6$, $E_7$ and $E_8$. \textit{Asymptotically locally flat} (ALF) metrics, on another hand, approach by definition the flat metric of an $S^1$-fibration over $\mathbb{R}^3$ or $\mathbb{R}^3/\mathbb{Z}_2$ minus a large geodesic ball, with fibers of asymptotically constant length. The boundary of the geodesic ball is diffeomorphic to $S^3/\Gamma$, where $\Gamma$ is either the the cyclic subgroup of $SU(2)$ in the first case, or the binary dihedral one in the second (see \textit{e.g.}~\cite{MR2778451}). 

In \cite{MR992334} Kronheimer gave a full classification of ALE gravitational instantons, constructing them by means of the hyperk\"ahler quotient construction of \cite{MR877637} and showing, moreover, that they are diffeomorphic to minimal resolutions of the Kleinian singularities $\mathbb{C}^2/\Gamma$. Although otherwise very powerful, the hyperk\"ahler quotient approach is not very good at producing explicit metrics. Such metrics have nevertheless been constructed in the $A_k$ case earlier on by Hitchin in \cite{MR520463}, who used a different, twistor-theoretic route, on minimal resolutions of Kleinian singularities of corresponding type. Despite the fact that the metrics he derived in this way, the so-called multi-center Eguchi-Hanson metrics, were already known from the work of Gibbons and Hawking \cite{Gibbons:1979zt}, Hitchin's approach shone a powerful conceptual light on the problem and contained in it the seeds of Kronheimer's generalization. The approach has more recently been extended in part to the $D_k$ case in \cite{MR2177322}, but the metric formulas derived there are somewhat unwieldy and still in a largely implicit form.

Geometrically, the principal difference between the $A_k$ and $D_k$ metrics lies in their symmetries: a tri-Hamiltonian one in the first case, respectively a hidden one of precisely the type we have considered in the second. So in the same way as the Gibbons-Hawking Ansatz can be used to construct explicitly the $A_k$ metrics, we use the quaternionic-based Ansatz we have devised to construct the $D_k$ ones. For that, we resort to the methods developed by Cherkis\pt-Kapustin \cite{MR1693628, MR1700937} and Cherkis\pt-Hitchin \cite{MR2177322} to convert the algebraic data of the universal deformation of the Klein quotient on which the gravitational instanton is defined into Lindstr\"om-Ro\v{c}ek twistor space data, which consists, essentially, of a single holomorphic transition function. Then, in combination with the Ansatz, the computation of the $D_k$ metrics reduces to the evaluation of several elliptic integrals, specifically one less than in \cite{MR2177322}\,---\,a fact which allows us to treat them in a uniform manner.



\section{Hyperk\"ahler spaces and twistors} 

The tangent bundle of a hyperk\"ahler manifold $M$ with standard triplet of symplectic 2-forms $\omega_1$, $\omega_2$, $\omega_3$ carries a representation of the algebra of imaginary quaternions generated by the tangent bundle endomorphisms $I_1 = \omega_3^{-1}\omega_2$, $I_2 = \omega_1^{-1}\omega_3$, $I_3 = \omega_2^{-1}\omega_1$. The closure property of the 2-forms ensures that these form integrable complex structures on $M$ \cite{MR887284}. In fact, any such manifold automatically possesses not just three but an entire $S^2$-family of integrable complex structures corresponding to all real-linear combinations of the form $x_1I_1+x_2I_2+x_3I_3$ with $x_1^2+x_2^2+x_3^2=1$. 

These can be assembled into a single integrable complex structure on the direct product space $Z = M \times S^2$ with the resulting complex manifold being known as the \textit{twistor space} of $M$ \cite{MR664330,MR877637}. By virtue of the isomorphism $S^2 \cong \mathbb{CP}^1$ the 2-sphere may be endowed with a complex structure of its own, and the natural projection $\pi: Z \rightarrow \mathbb{CP}^1$ is then holomorphic. As  a complex manifold, $\mathbb{CP}^1$ can be obtained by patching together two copies of the complex plane $\mathbb{C}$ with coordinates $\zeta$ and $\tilde{\zeta}$ related on the intersection domain by the holomorphic transition relation $\tilde{\zeta} = 1/\zeta$. In these notes we shall choose the coordinate charts in such a way that $\zeta = 0$ and $\infty$ correspond through the stereographic map realizing explicitly the isomorphism to the ``north" and ``south" poles on the sphere, that is, the points $(x_1,x_2,x_3) = (0,0,1)$ and $(0,0,-1)$, respectively\,---\,or, equivalently put, to the conjugate complex structures $I_3$ and $-I_3$ on the sphere of complex structures. 

Besides the complex structure, the twistor space $Z$ is also endowed, very importantly, with a \textit{real structure}\,---\,an anti-holomorphic involution induced by antipodal conjugation on $S^2$ taking $(p,\zeta) \mapsto (p,-1/\bar{\zeta})$ for any point $p \in M$. In the main body of the paper we will often use for antipodally-conjugate variables the superscript notation $\zeta^c \equiv -1/\bar{\zeta}$. 

Observe also that the above definition of hyperk\"ahler manifolds implies that they are naturally Riemannian, or perhaps pseudo-Riemannian, with a metric given by $g(X,Y) = -\allowbreak\omega_1(X,I_1Y) = -\allowbreak\omega_2(X,I_2Y) =  -\allowbreak\omega_3(X,I_3Y)$ for any pair of vector fields $X,Y \in TM$. This hyperk\"ahler metric together with the Fubini-Study metric on $S^2$ induce in turn a metric on the twistor space. When coupled to the twistor complex structure the twistor metric yields a natural 2-form on $Z$. This form is not of K\"ahler type in general, rather the twistor space metric is \textit{balanced} \cite{MR623721,MR1669956}. This fact, however, will not play a significant role in our considerations here.

\section{The generalized Legendre transform construction} 

\subsection{}

The class of hyperk\"ahler spaces that we shall be concerned with in this work has the distinguishing feature that the holomorphic projection $\pi$ factorizes completely through a direct sum of line bundles over $\mathbb{CP}^1$ of positive even degree:
\begin{center}
\begin{tikzcd}
Z \rar & \displaystyle \bigoplus_{\I=1}^{\dim_{\pt \mathbb{H}} M} \! \mathcal{O}(2j_{\I}) \rar & \mathbb{CP}^1
\end{tikzcd}
\end{center}
In terms of the associated symmetry properties of the hyperk\"ahler metric there is an important qualitative difference between the intermediate summands with $j_{\I} =1$ and the ones with $j_{\I} \geq 2$: each summand of the first type corresponds to a tri-Hamiltonian\footnote{A tri-Hamiltonian action on a hyperk\"ahler space is by definition an action which is simultaneously Hamiltonian with respect to each of the three fundamental hyperk\"ahler symplectic forms $\omega_1$, $\omega_2$, $\omega_3$.} (and hence Killing) \textit{vector} symmetry, while each summand of the second type is associated to a certain kind of Killing \textit{tensor} symmetry (or Killing \textit{spinor} symmetry, in four dimensions). The approach which produces such metrics was discovered by Lindstr\"om and Ro\v{c}ek and is known as the Legendre transform construction when all the intermediate line bundles have degree 2 \cite{MR710273}, and as the \textit{generalized} Legendre transform construction when at least one of them has degree greater than or equal to 4 \cite{MR929144} (for a unified and in fact even more generic perspective see also the more recent reference \cite{Lindstrom:2008gs}). Both constructions have emerged originally in connection to supersymmetric quantum field-theoretic problems and were later on given mathematical twistor-theoretic interpretations in \cite{MR877637} and \cite{MR1447294}, respectively. The first and simplest one yields the well-known class of toric hyperk\"ahler metrics with a free local tri-Hamiltonian $\mathbb{R}^n$-action of rank $n$ equal to the quaternionic dimension of the space, $\dim_{\pt \mathbb{H}} M$, described by the Gibbons-Hawking Ansatz \cite{Gibbons:1979zt} or one of its higher-di\-men\-sio\-nal analogues \cite{MR953820}, depending on whether $n=1$ or $>1$. 

We begin our investigations with a detailed unified review of these approaches aimed at establishing a foundation for our subsequent considerations. For simplicity we assume that all the positive integers $j_{\I}$ take the same value, which we denote by $j$, and although we differentiate at times between the case when $j=1$ and the case when $j \geq 2$, we treat them together inasmuch as possible. The generic case involves only a simple extension of the discussion which follows, entailing essentially a further refinement of the indices to accommodate for mixed-type combinations. 

Given a hyperk\"ahler space $M$ with twistor space $Z$, there exists a natural global holomorphic section of the bundle $\Lambda^2T_F^* \otimes \pi^*\mathcal{O}(2)$ over $Z$, where $T_F = \ker d\pi$ is the holomorphic tangent bundle along the fibers, $T^*_F$ its dual, and $\pi^*\mathcal{O}(2)$ is the pullback over $Z$ of the $\mathcal{O}(2)$ bundle over $\mathbb{CP}^1$. Let $N$ and $S$ be two open sets in $Z$ projecting down to the standard two-set open covering of $\mathbb{CP}^1$.  In the local trivialization over $N$ this section takes the form
$\omega_N(\zeta) = \omega_+ + \omega_{\hp 0} \hp \zeta + \omega_-\zeta^2$,
where $\omega_{\pm} = \pm \frac{1}{2}(\omega_1 \pm i \hp \omega_2)$ and $\omega_{\hp 0} = \omega_3$ form a frame in the complexified linear space spanned by the fundamental hyperk\"ahler symplectic forms. Note that thus defined the elements of this frame satisfy the \textit{alternating conjugation property} $\bar{\omega}_m = (-)^m \omega_{-m}$, where $m = -1,0,+1$. One has also a corresponding component of the section over $S$, and on the overlap $N \cap S$ the  two components are patched together by means of the usual $\mathcal{O}(2)$ transition function, $\omega_S(\tilde{\zeta}) = \zeta^{-2} \omega_N(\zeta)$. The alternating conjugation property of its coefficients translates into the reality property $\overline{\omega_N(\zeta^c)} = - \pt \omega_S(\tilde{\zeta})$ for the section in which the two components are interchanged. 

For each $\zeta$ and $\tilde{\zeta}$, $\omega_N(\zeta)$ and respectively $\omega_S(\tilde{\zeta})$ define complex symplectic forms on the corresponding twistor fibers. Thus, the twistor fiber over a generic point in $\mathbb{CP}^1$ can be viewed as a copy of the manifold $M$ endowed with the particular hyperk\"ahler complex structure labeled by that point, as well as with a complex symplectic form holomorphic with respect to it. By the complex version of Darboux's theorem we can then define for each fiber's symplectic form a local Darboux coordinate chart, $\omega_N(\zeta) = dp_{\I,N}(\zeta) \wedge dx^{\I}_N(\zeta)$ and $\omega_S(\tilde{\zeta}) = dp_{\I,S}(\tilde{\zeta}) \wedge dx^{\I}_{S}(\tilde{\zeta})$, where summation over the repeated index $I$ is implied. The symplectic coordinates depend holomorphically on the local $\mathbb{CP}^1$ coordinate, and can in fact be seen as local holomorphic coordinates on the twistor space. Over the intersection domain their symplectic gluing may be described by a (twisted\,---\,to account for the $\mathcal{O}(2)$ twist of the holomorphic twistor space 2-form) canonical transformation, with a generating function of, say, type II: 
{\allowdisplaybreaks
\begin{equation}
\begin{aligned}
p_{\I,N} & = \zeta \frac{\partial F_2(x_N,p_S,\zeta)}{\partial x^{\I}_N} \\[-1pt]
x^{\I}_S & = \frac{1}{\zeta} \frac{\partial F_2(x_N,p_S,\zeta)}{\partial p_{\I,S}} \rlap{.}
\end{aligned}
\end{equation}
}%
In the generalized Legendre transform construction one assumes that two coordinate charts and a single patching function between them suffice, although one typically allows the coordinates to be possibly multi-valued and to have singularities or branching points\,---\,except, crucially, at points in $Z$ above $\zeta = 0$ and $\tilde{\zeta} =0$. Another important assumption is that this generating function is of the form 
\begin{equation}
F_2(x_N,p_S,\zeta) = \zeta^{1-2j} p^{\phantom{\I}}_{\I,S} \hp x^{\I}_N + i \hp H(x_N/\zeta^j,\zeta)
\end{equation}
where $H$ is a holomorphic function of its variables, which do not include the $p_{\I,S}$, with possible singularities and branching points. For this choice the gluing equations read 
\begin{align}
\zeta^j x^{\I}_S & = \zeta^{-j} x^{\I}_N \equiv x^{\I}_{NS} \label{arct-trop-anta} \\
\zeta^{1-j}p_{\I,S} & = \zeta^{j-1} p_{\I,N} - i \hp \frac{\partial H(x_{NS},\zeta)}{\partial x^{\I}_{NS}} \rlap{.} \label{pN-pS}
\end{align}
The first equation asserts essentially that for each value of the index $I$, $x^{\I}_N$ and $x^{\I}_S$ represent the \textit{arctic} and \textit{antarctic} components of a section of the $\pi^*\mathcal{O}(2j)$ bundle over $Z$. It is useful to define also a \textit{tropical} component, $x^{\I}_{NS}$, supported on $N \cap S$. The same geographic monikers will be used throughout the paper to refer to similar components of sections of line bundles of even degree over $\mathbb{CP}^1$. 

Finally, consistently with the reality properties of the twistor space holomorphic 2-form, we assume furthermore that the coordinates satisfy well-defined reality properties of their own: 
\begin{equation} \label{real-xp}
\begin{aligned}
\overline{x^{\I}_N(\zeta^c)} & = (-)^j \, x^{\I}_S(\tilde{\zeta}) \\[2pt]
\overline{p_{\I,N}(\zeta^c)} & = (-)^{1-j} \, p_{\I,S}(\tilde{\zeta}) \rlap{.}
\end{aligned}
\end{equation}
In this paper we will refer to the substitution $\zeta \mapsto \zeta^c$ followed by complex conjugation as \textit{antipodal conjugation}. So, in this language, we assume that the $x^{\I}$ and $p_{\I}$ twistor coordinates are, up to a possible sign, antipodally self-conjugate.

From the regularity conditions at $\zeta=0$ and the first gluing equation we have for the antarctic components Taylor $\zeta$-expansions of the form
{\allowdisplaybreaks
\begin{align}
x^{\I}_N(\zeta) & = \sum_{n=-j}^{j} x^{\I}_{n} \hp \zeta^{j-n} = x^{\I}_j + \mathcal{O}(\zeta)  \\[-2pt]
p_{\I,N}(\zeta) & = \sum_{n=-\infty}^{1-j} p_{\I,n} \zeta^{1-j-n} = p_{\I,1-j} + \mathcal{O}(\zeta) \rlap{.}
\end{align}
}%
Note that the leading terms in the expansions give a complete set of holomorphic coordinates with respect to the complex structure labeled on the twistor sphere by $\zeta=0$, that is, $I_3$, and in what follows we sometimes denote them alternatively by $x^{\I}_j \equiv z^{\I}$ and $p_{\I,1-j} \equiv u_{\I}$ to highlight this special role. 

The next step is to equate the coefficients of the powers of $\zeta$ on the two sides of the gluing equations. For the first equation, due to the finiteness of the Taylor expansion this is a straightforward proposition. In the second equation one can accomplish this with the help of contour integration. After multiplying the equation with $\zeta^{m-1}$ we integrate the result over a contour contained initially in the tropical region of the Riemann sphere. By appropriately deforming the contour for each term in part, the integral will pick the corresponding coefficients in the Laurent expansions of the first two terms while providing a contour-integral prescription for the contribution of the third term, whose holomorphic structure may vary from case to case. The outcome of these manipulations can be expressed in terms of the single real-valued function
\begin{equation} \label{L-oint}
L = \oint_{\mathcal{C}} \frac{d\zeta}{2\pi i \hp \zeta} H(x_{NS}(\zeta),\zeta)
\end{equation}
depending on the coefficients $x^{\I}_n$. The reality conditions satisfied by the $x\hp$-variables ensure that we can always choose the contour $\mathcal{C}$ in such a way that the result of the integration is real. In the end, the matching of powers of $\zeta$ in the gluing equations gives, respectively,
\begin{align}
(-)^m \bar{x}^{\I}_{-m} & = x^{\I}_m \label{alt-conj-x} \\[-2pt]
(-)^m \bar{p}_{\I,-m} & = p_{\I,m} - i \frac{\partial L}{\partial x^{\I}_{-\rlap{$\scriptstyle m$}}} 
\end{align}
where $m = -j, \dots, j$. The first set of equations is simply a reflection of the reality properties of the $\mathcal{O}(2j)$ sections. The second set of equations is the one storing the interesting geometric information of the manifold, which we can now see is encoded in the holomorphic structure of the function $H$ through the agency of the real-valued contour integral $L$. Taking a closer look, we notice that for any value of $j$ the $m=-j$ component reads
\begin{equation}
\frac{\partial L}{\partial x^{\I}_j} = - i \hp p_{\I,-j}
\end{equation}
whereas for the remaining components we have
{\allowdisplaybreaks
\begin{align}
& \text{for $j=1$\pt:} \qquad \frac{\partial L}{\partial x^{\I}_0} = 2 \hp \Im \hp u_{\I} \label{LT-eqs} \\
& \text{for $j \geq 2$\pt:} \qquad \hp \frac{\partial L}{\partial x^{\I}_{\rlap{$\scriptstyle m$}}} = \pt
\begin{cases}
- i \hp u_{\I} & \text{for $m=j-1$} \\[2pt]
\hspace{10pt} 0 & \text{for $m=j-2,\dots,0$} \rlap{.}
\end{cases} \label{GLT-eqs}
\end{align}
}%
The equations corresponding to negative values of $m$ follow trivially from these through complex conjugation.

We can have a deeper understanding of these differential equations if we introduce the following real-valued function
{\allowdisplaybreaks
\begin{align}
& \text{for $j=1$\pt:} \qquad \kappa(z,u,\bar{z},\bar{u}) = \langle \pt  L - 2 \hp \Im \hp u_{\I} \pt x^{\I}_0 \pt \rangle_{x^{\I}_0} \\
& \text{for $j\geq 2$\pt:} \qquad \kappa(z,u,\bar{z},\bar{u}) = \langle \pt  L  + i \hp u_{\I}x^{\I}_{j-1} - i \hp \bar{u}_{\I} \bar{x}^{\I}_{j-1} \pt \rangle_{x^{\I}_{j-1},.., \pt x^{\I}_0,.., \pt x^{\I}_{1-j}}
\end{align}
}%
where the chevron brackets signify that we perform a simultaneous extremization with respect to the variables listed in the lower right-hand side subscript. When multiple extrema are possible we assume a selection criterion exists rendering the result unique. For $j=1$ this is a regular multivariate Legendre transform, and the $j\geq 2$ cases can be viewed as generalizations. Observe that the equations \eqref{LT-eqs} and \eqref{GLT-eqs} represent precisely the extremizing conditions associated to these transforms which allow a trade-off of variables. One can in principle solve them to express the coefficients $x^{\I}_m$ as implicit functions of the complex variables $z^{\I}$, $u_{\I}$ (and their complex conjugates). Recalling that the latter are holomorphic with respect to the complex structure $I_3$, the claim is then that the resulting function $\kappa(z,u,\bar{z},\bar{u})$  provides a K\"ahler potential for the corresponding K\"ahler form $\omega_3$. 

To see that, observe that in either case above one has
\begin{equation} \label{delK}
\frac{\partial \kappa}{\partial z^{\I}}  = \frac{\partial L}{\partial x^{\I}_{j}}  = - i \hp p_{\I,-j}
\qquad\text{and}\qquad
\frac{\partial \kappa}{\partial u_{\I}} \hspace{5pt} = i \hp x^{\I}_{j-1} \rlap{.}
\end{equation}
On the other hand, by Taylor expanding in $\zeta$ the Darboux formula for $\omega_N(\zeta)$ we can identify 
$\omega_+ = dp^{\phantom{\I}}_{\I,1-j} \wedge dx^{\I}_j$ and $\omega_3 = dp^{\phantom{\I}}_{\I,-j} \wedge dx^{\I}_{j} + dp^{\phantom{\I}}_{\I,1-j} \wedge dx^{\I}_{j-1}$. 
In terms of the $I_3$-holomorphic coordinates the $(2,0)$-form can be transcribed trivially to
\begin{equation}
\omega_+ = du_{\I} \wedge dz^{\I}
\end{equation}
while using the previous observation for the $(1,1)$ form we can write successively
\begin{equation}
\omega_3 = d(p_{\I,-j}dz^{\I} - x^{\I}_{j-1} du_{\I}) = i \hp d\bigg(\frac{\partial \kappa}{\partial z^{\I}} dz^{\I} + \frac{\partial \kappa}{\partial u_{\I}} du_{\I} \! \bigg) = - i \hp \partial_{\scriptscriptstyle I_3}\bar{\partial}_{\scriptscriptstyle I_3} \kappa \rlap{,}
\end{equation}
which demonstrates the claim.

\subsection{}

Possession of a K\"ahler potential and of an accompanying set of complex coordinates allows us to derive formulas for the metric and corresponding K\"ahler form. Before we proceed to present the details of this calculation let us introduce a number of notations. First, we will denote derivatives with respect to the coefficients $x^{\I}_m$ by means of indices. For example, the Hessian matrix of $L$ will be denoted by $(L_{x^{\I}_{m}x^{\J}_{\smash{m'}}})^{\phantom{I}}_{m,m'=-j,\dots,j}$. Secondly, we will introduce special notations for some of these coefficients to indicate that we consider only a limited range of them in terms of the values that the integer $m$ can take. These notations are summed up in Table \ref{Indices}.
\begin{table}[ht]
\centering
\begin{tabular}{l|l}
Symbol & Range \\[2pt] \hline 
$a$, $b$ & $(x^{\I}_m)^{\phantom{I}}_{m=j-1,\dots,1-j}$ \rule{0pt}{2.5ex} \\[2pt]
$a'$ & The range of $a$ without $x^{\I}_{j-1}$ \\[2pt]
$b''$ & The range of $\hp b\hp$ without $x^{\I}_{1-j}$ \\[2pt] \hline
\end{tabular}
\caption{\protect\rule{0pt}{2.5ex}} \label{Indices}
\vspace*{-12pt}
\end{table}
In the course of the calculation we will also need to consider certain square matrices issuing through cuts and inversions from the Hessian matrix of $L$. To ease references to them we group their definitions in Table \ref{Submatrices}.
\begin{table}[ht]
\centering
\begin{tabular}{l@{}l|l}
\multicolumn{2}{l|}{Matrix definition}   & \,Inverse matrix \\[2pt] \hline
$L_{ab}$ 						& $\coloneqq (L_{x^{\I}_{m}x^{\J}_{\smash{m'}}})^{\phantom{I}}_{m,m'=j-1,\dots,1-j}$ 	& \,$L^{ab}$ 	     \rule{0pt}{2.8ex}	\\[3pt]
$\Phi^{\scriptscriptstyle I\bar{J}}$ 	& $\coloneqq$  the submatrix $L^{x^{\I}_{j-1}x^{\J}_{1-j}}$ of $L^{ab}$ 			& \,$\Phi_{\scriptscriptstyle \bar{J}I}$ \\[4pt]
$M_{a'b''}$ 					& $\coloneqq$ the submatrix $L_{a'b''}$ of $L_{ab}$ 							& \,$M^{b''a'}$ 					\\[4pt]
$N_{b''a'}$	 				& $\coloneqq$ the submatrix $L_{b''a'}$ of $L_{ab}$ 							& \,$N^{a'b''}$ 					\\[4pt] \hline
\end{tabular}
\caption{\protect\rule{0pt}{2.5ex}} \label{Submatrices}
\vspace*{-12pt}
\end{table}

From either the equations \eqref{LT-eqs} or equations \eqref{GLT-eqs}, with the notations outlined in the first rows of both Table \ref{Indices} and Table \ref{Submatrices}, for any value of $j$, including $j=1$, we have
\begin{equation}
\begin{aligned}
\frac{\partial a}{\partial z^{\I}} & = - L^{ab} L_{b \hp x^{\I}_j} & \frac{\partial a}{\partial u_{\I}} & = - i \hp L^{ax^{\I}_{j-1}} \\
\frac{\partial a}{\partial \bar{z}^{\I}} & = (-)^{j-1} L^{ab} L_{b \hp x^{\I}_{-j}} \quad & \frac{\partial a}{\partial \bar{u}_{\I}} & = (-)^{j-1} i \hp L^{a x^{\I}_{1-j}} \rlap{.}
\end{aligned}
\end{equation}
Acting then with antiholomorphic derivatives on the first derivative expressions \eqref{delK} for the K\"ahler potential and using these formulas either directly or in combination with the chain rule we obtain
\begin{align}
\left(
\begin{array}{c|c}
\kappa_{z^{\I}\bar{z}^{\J}} & \kappa_{z^{\I}\bar{u}_{\J}} \!\! \\[5pt] \hline
\kappa_{u_{\I}\bar{z}^{\J}} & \kappa_{u_{\I}\bar{u}_{\J}} \rule{0pt}{12pt} \!\!
\end{array}
\right)
=  (-)^j \!
\left(
\begin{array}{c|c}
\!\! L_{x^{\I}_jx^{\J}_{-j}} \! - L_{x^{\I}_ja} L^{ab} L_{b \hp x^{\J}_{-j}} \! & -i \hp L_{x^{\I}_ja} L^{ax^{\J}_{1-j}} \!\! \\[6pt] \hline
- i \hp L^{x^{\I}_{j-1}b} L_{b \hp x^{\J}_{-j}} & L^{x^{\I}_{j-1}x^{\J}_{1-j}} \rule{0pt}{15pt} 
\end{array}
\right) \rlap{.}
\end{align}

It turns out that the $\kappa_{z^{\I}\bar{z}^{\J}}$ component can be cast in a more convenient form which has the effect of making the hyperk\"ahler structure of the metric more readily transparent. To write down this form let us define a matrix $\Phi^{\scriptscriptstyle I\bar{J}}$ as in the second row of Table \ref{Submatrices}.  Notice that we have just shown that $\Phi^{\scriptscriptstyle I\bar{J}} = (-)^j\kappa_{u_{\I}\bar{u}_{\J}}$ and so clearly, since the K\"ahler potential is real, this matrix is Hermitian. Then so must be its inverse, $\Phi_{\scriptscriptstyle \bar{J}I}$. Let $\Phi_{\scriptscriptstyle J\bar{I}}$ denote the elements of the latter's complex conjugate matrix. With these definitions in place we claim then that
\begin{equation} \label{zzbar-term}
L_{x^{\I}_jx^{\J}_{-j}} \! - L_{x^{\I}_ja} L^{ab} L_{b \hp x^{\J}_{-j}} = \Phi_{\scriptscriptstyle I\bar{J}} - L_{x^{\I}_ja} L^{a \hp x^{P}_{1-j}} \Phi_{\scriptscriptstyle \bar{P}Q} L^{x^{Q}_{j-1} b} L_{b \hp x^{\J}_{-j}} \rlap{.}
\end{equation}

To prove the claim we begin by noticing that the components of the matrices $L^{ab}$ and $L^{a \hp x^{P}_{1-j}} \Phi_{\scriptscriptstyle \bar{P}Q} L^{x^{Q}_{j-1} b}$ with $a = x^{\I}_{j-1}$ and any value of $b$ are identical, and likewise those with $b = x^{\I}_{1-j}$ and any value of $a$, so the corresponding terms from the two sides of equation \eqref{zzbar-term} can be dropped out to yield the equivalent statement
\begin{equation} \label{intermed-claim}
L_{x^{\I}_jx^{\J}_{-j}} \! - L_{x^{\I}_ja'} L^{a'b''} L_{b'' x^{\J}_{-j}} = \Phi_{\scriptscriptstyle I\bar{J}} - L_{x^{\I}_ja'} L^{a' x^{P}_{1-j}} \Phi_{\scriptscriptstyle \bar{P}Q} L^{x^{Q}_{j-1} b''} L_{b'' x^{\J}_{-j}} 
\end{equation}
where we have indicated restrictions in the ranges of indices by means of the notations defined in the last two rows of Table \ref{Indices}. Consider now the two square submatrices of $L_{ab}$ defined in the last two rows of Table \ref{Submatrices}. Their matrix inverses satisfy 
\begin{equation}
\begin{aligned}
N^{a'b''} & = L^{a'b''} \! - L^{a' x^{P}_{1-j}} \Phi_{\scriptscriptstyle \bar{P}Q} L^{x^{Q}_{j-1} b''} \\
\Phi_{\scriptscriptstyle I\bar{J}} & = L_{x^{\I}_{j-1}x^{\J}_{1-j}} \! - L_{x^{\I}_{j-1}b''} M^{b''a'} L_{a' x^{\J}_{1-j}} \rlap{.}
\end{aligned}
\end{equation}
Using these identities we can then see that the equation \eqref{intermed-claim} is in turn equivalent to 
\begin{equation} \label{twosubms}
L_{x^{\I}_jx^{\J}_{-j}} \! - L_{x^{\I}_ja'} N^{a'b''} L_{b'' x^{\J}_{-j}} = L_{x^{\I}_{j-1}x^{\J}_{1-j}} \! - L_{x^{\I}_{j-1}b''} M^{b''a'} L_{a' x^{\J}_{1-j}} \rlap{.}
\end{equation}
Observe that (see \textit{e.g.} the formulas in part (ii) of Theorem 2.2 in \cite{MR1873248}) the two sides of this equation come up in the block-matrix inversions
{\allowdisplaybreaks
\begin{equation} \label{bl-diag-inv}
\begin{aligned}
\left(
\begin{array}{c|c}
\hspace{5.5pt} L_{x^{\I}_ja'} \hspace{5.5pt} & \hspace{6.5pt} L_{x^{\I}_jx^{\J}_{-j}} \hspace{6.5pt} \!\! \\[5pt] \hline
L_{b''a'} & L_{b'' x^{\J}_{-j}} \rule{0pt}{12pt} \!\!
\end{array}
\right)^{\!\!-1}
& \! =
\left(
\begin{array}{c|c}
\cdots & \cdots \\[5pt] \hline
(\text{l.h.s.\ of eq.\,\eqref{twosubms}})^{-1} \hspace{1.2pt} & \! \cdots \rule{0pt}{12pt} \!\!
\end{array}
\right) \\
\left(
\begin{array}{c|c}
L_{x^{\I}_{j-1}b''} & L_{x^{\I}_{j-1}x^{\J}_{1-j}} \!\! \\[5pt] \hline
L_{a'b''} & L_{a' x^{\J}_{1-j}} \rule{0pt}{12pt} \!\!
\end{array}
\right)^{\!\!-1}
& \! =
\left(
\begin{array}{c|c}
\cdots & \cdots \\[5pt] \hline
(\text{r.h.s.\ of eq.\,\eqref{twosubms}})^{-1} & \! \cdots \rule{0pt}{12pt} \!\!
\end{array}
\right) 
\end{aligned}
\end{equation}
}%
Recalling the definitions of the primed indices one can easily see that the matrices being inverted are
\begin{equation} \label{two-submat}
(L_{x^{\I}_{m}x^{\J}_{\smash{m'}}})_{\hp \begin{subarray}{l} m\phantom{'} = j, \dots, 2-j \\ m' = j-2, \dots, -j \end{subarray}} 
\qquad\text{and}\qquad
(L_{x^{\I}_{m}x^{\J}_{\smash{m'}}})_{\hp \begin{subarray}{l} m\phantom{'} = j-1, \dots, 1-j \\ m' = j-1, \dots, 1-j \end{subarray}}
\end{equation}
respectively. 

So far all we have done was to perform a series of rather dry and technical linear algebra manipulations. Now comes the most important part of the argument, whose relevance transcends the limited scope of this proof and represents a well-known signature of the twistor approach. The contour-integral expression \eqref{L-oint} implies that the function $L$ satisfies the following second-order differential equations
\begin{equation} \label{Bogo-equiv-j=2}
L_{x^{\I}_mx^{\J}_{\smash{m'}}} = L_{x^{\I}_{\smash{m+k}}x^{\J}_{\smash{m'\!-k}}} \! =  L_{x^{\J}_{m}x^{\I}_{\smash{m'}}}
\end{equation}
for all admissible values of $m$, $m'$ and $k$. This is the concrete expression of the celebrated twistor correspondence between differential equations and holomorphicity. In particular, it allows us to define the functions
\begin{equation} \label{h-bar-def}
\bar{h}_{\I\J,n} = L_{x^{\I}_{k}x^{\J}_{n-k}}
\end{equation}
independent of $k$, with the complex conjugation having been chosen for ulterior convenience. The reality properties \eqref{alt-conj-x} of the $x^{\I}_m$ variables together with the reality of $L$ imply the alternating conjugation property $\bar{h}_{\I\J,n} = (-)^n h_{\I\J,-n}$. Also, from \eqref{L-oint} one has the following contour-integral representation:
\begin{equation} \label{h-bar-oint}
\bar{h}_{\I\J,n} = \oint_{\mathcal{C}} \frac{d\zeta}{2\pi i} \, \zeta^{n-1} \frac{\partial^2 H(x_{NS},\zeta)}{\partial x^{\I}_{NS}\partial x^{\J}_{NS}} \rlap{.}
\end{equation}

Notice that because of these identifications the two matrices \eqref{two-submat} are one and the same and equal to $(\hp \bar{h}_{\I\J, \hp m+m'} \hp )_{\pt 1-j \leq m,m' \leq j-1}$. But in the light of the formulas \eqref{bl-diag-inv} this is a sufficient condition for the equation \eqref{twosubms} to hold and thus, following backwards the chain of equivalences, the claim is proved. 

Assembling everything together, the result of the calculations can then be described as follows: 

\begin{proposition} \label{GLT-formulas}
The (generalized) Legendre transform construction gives the following set of formulas for the hyperk\"ahler metric
\begin{equation}
ds^2 = (-)^j [ \pt dz^{\I} \Phi_{\scriptscriptstyle I\bar{J}} \pt d\bar{z}^{\scriptscriptstyle \bar{J}} + (du_{\I} + dz^{\scriptscriptstyle P} \! A_{\scriptscriptstyle PI}) \Phi^{\scriptscriptstyle I\bar{J}} (d\bar{u}_{\scriptscriptstyle \bar{J}} + \bar{A}_{\scriptscriptstyle \bar{J}\bar{Q}} d\bar{z}^{\scriptscriptstyle \bar{Q}}) \pt ]
\end{equation}
and hyperk\"ahler symplectic forms
\begin{gather}
\omega_+ = du_{\I} \wedge dz^{\I} \\
\omega_3 = (-)^{j-1} i \pt [ \pt \Phi_{\scriptscriptstyle I\bar{J}} \hp dz^{\I} \!\wedge d\bar{z}^{\scriptscriptstyle \bar{J}} + \Phi^{\scriptscriptstyle I\bar{J}} (du_{\I} + dz^{\scriptscriptstyle P} \! A_{\scriptscriptstyle PI}) \hspace{-1pt} \wedge \hspace{-1pt} (d\bar{u}_{\scriptscriptstyle \bar{J}} + \bar{A}_{\scriptscriptstyle \bar{J}\bar{Q}} d\bar{z}^{\scriptscriptstyle \bar{Q}}) \pt ] \rlap{,}
\end{gather}
where the Hermitian $\Phi$-matrix and its inverse are defined in the second row of Table \ref{Submatrices}, and 
$A_{\I\J} = - i \hp L_{x^{\I}_ja} L^{a \hp x^{P}_{1-j}} \Phi_{\scriptscriptstyle \bar{P}J}$,
with $\bar{A}_{\scriptscriptstyle \bar{J}\bar{I}} = -i \pt \Phi_{\scriptscriptstyle \bar{J}P} L^{x^P_{j-1}a}L_{a x^{\I}_{-j}}$ its complex conjugate.
\end{proposition}

Let us make a few remarks in connection to this Proposition.
\begin{itemize}
\setlength\itemsep{3pt}
\item[1.] The formulas are valid for both $j=1$ and $j \geq 2$. 
\item[2.] They are expressed in a coordinate frame holomorphic with respect to a \textit{particular} hyperk\"ahler complex structure out of the two-sphere's worth of these available, in this case $I_3$. Nevertheless, note that this does not come at the expense of obscuring the hyperk\"ahler quaternionic structure which is actually quite transparent.
\item[3.] In the generalized Legendre transform construction the hyperk\"ahler structure is completely determined by a subset of the second derivatives of a single real-valued function $L$, more specifically by the functions $(h_{\I\J,n})_{n=1-2j,\dots,2j-2}$. As we will see later on, a detailed analysis of the four-dimensional case with $j=2$ strongly points to a further reduction in the number of these needed for a complete description, the key to which involves choosing a different coordinate system. We conjecture that the minimum subset necessary is given by $(h_{\I\J,n})_{n=-j,\dots,j}$. 
\item[4.] However, observe that the potential function $L$ depends on the variables $x^{\I}_m$ while the metric is given in a coordinate frame determined by the very different variables $z^{\I}$ and $u_{\I}$. The connection between these two sets of variables is provided by the equations \eqref{LT-eqs}--\eqref{GLT-eqs} and is in most cases implicit rather than explicit. Apart from the exceptional case $j=1$, this poses very serious practical difficulties in the way of writing down explicit metrics. We will see that the different coordinate system hinted to above brings substantial simplifications in this respect as well. 
\end{itemize}

\subsection{}

The metrics in this class possess automatically certain symmetries, of either classical or hidden type, depending on whether $j=1$ or $j \geq 2$, respectively. 

For $j=1$ the metrics have a well-known abelian tri-Hamiltonian\,---\,and hence Killing\,---\,symmetry of rank equal to the quaternionic dimension of the space, generated by the vector fields
\begin{equation}
X^{\I}_0 = \frac{\partial}{\partial u_{\I}} + \frac{\partial}{\partial \bar{u}_{\I}} \rlap{,}
\end{equation}
with the role of Hamiltonian functions being played by the components of the $\mathcal{O}(2)$ sections $x^{\I}$. That is, we have
\begin{equation}
\iota_{X^{\I}_0} \omega^{\phantom{\I}}_m = dx^{\I}_m
\end{equation}
for $m=-1,0,+1$, and so the Lie action of each $X^{\I}_0$ preserves the three hyperk\"ahler symplectic forms separately:
\begin{equation}
\mathcal{L}_{X^{\I}_0} \omega_m = 0 \rlap{.}
\end{equation}

The metrics with $j \geq 2$ are characterized by a \textit{twisted} or $\textit{hidden}$ version of this symmetry. To see that, let us consider the vector fields $\mathcal{X}^{\I}_m$ with $m = 0,\dots,j-1$, of type $(1,0)$ relative to the complex structure $I_3$ with respect to which the coordinates $u_{\I}, z^{\I}$ are holomorphic, and such that
\begin{equation}
\begin{aligned}
\iota_{\mathcal{X}^{\I}_m} \omega_+ & = \partial_{I_3} x^{\I}_{m+1} \\
\iota_{\mathcal{X}^{\I}_m} \omega_{\hp 0} \pt\hp & = \bar{\partial}_{I_3} x^{\I}_{m} 
\end{aligned}
\end{equation}
where $\partial_{I_3}$ denotes the Dolbeault operator corresponding to $I_3$. The first condition defines the vector fields as symplectic gradients, and the second one can be shown with some effort to follow from the differential constraints \eqref{Bogo-equiv-j=2}. Alternatively, these conditions can be understood as consequences of the holomorphicity of the $\mathcal{O}(2j)$ sections $x^{\I}$, in a twistor space sense. Let us assemble further the vector fields
\begin{equation}
X^{\I}_m = 
\begin{cases}
\mathcal{X}^{\I}_m & \text{for $0 < m \leq j-1$} \\[2pt]
\mathcal{X}^{\I}_0 + \bar{\mathcal{X}}^{\I}_0 & \text{for $m=0$} \\[2pt]
(-)^m \bar{\mathcal{X}}^{\I}_{-m} & \text{for $0 > m \geq 1-j$} \\[1.5pt]
0 & \text{otherwise}
\end{cases}
\end{equation}
with $2j-1$ non-vanishing components and satisfying the reality condition $\bar{X}^{\I}_m = (-)^m X^{\I}_{-m}$. Note, incidentally, that the first and last non-vanishing components are holomorphic respectively anti-holomorphic vector fields with respect to the complex structure $I_3$:
\begin{equation}
X^{\I}_{j-1} = \frac{\partial}{\partial u_{\I}},
\qquad \text{\dots} \qquad , \,
X^{\I}_{1-j} = (-)^{j-1} \frac{\partial}{\partial \bar{u}_{\I}} \rlap{.}
\end{equation}
Recalling that $\omega_+$, $\omega_{\hp 0}$, $\omega_-$ are of type $(2,0)$, $(1,1)$, $(0,2)$ relative to $I_3$, respectively, we infer then immediately that for all $m=-j,\dots,j$ we have
\begin{equation}
\iota_{X^{\I}_{m-1}} \omega_+ + \iota_{X^{\I}_{m\rule[0.8ex]{0pt}{0pt}}} \omega_{\hp 0} + \iota_{X^{\I}_{m+1}} \omega_- = dx^{\I}_m \rlap{.}
\end{equation}
Thus, the components of the sections $x^{\I}$ can be seen in this case as Hamiltonian functions, too, but only in this generalized sense. Acting with a total derivative on this equation and using Cartan's formula and the closure of the hyperk\"ahler symplectic forms gives eventually
\begin{equation}
\mathcal{L}_{X^{\I}_{m-1}} \omega_+ + \mathcal{L}_{X^{\I}_{m\rule[0.8ex]{0pt}{0pt}}} \omega_{\hp 0} + \mathcal{L}_{X^{\I}_{m+1}} \omega_- = 0 \rlap{.}
\end{equation}
Actions of this type on a hyperk\"ahler space have been studied in \cite{MR1848654} by Bielawski under the name of \textit{twistor group actions}.

\subsection{}

In four dimensions (\textit{i.e.}~$\dim_{\mathbb{H}} M = 1$) the indices $I,J,\dots$ drop out as they take only a single value and the formulas of Proposition \ref{GLT-formulas} read simply
\begin{equation} \label{dimH=1-GLT-metric}
ds^2 = (-)^j [\pt \Phi \hp |dz|^2 + \Phi^{-1}|du + A \hp dz|^2]
\end{equation}
respectively
\begin{gather} 
\omega_+ = du \wedge dz \\
\omega_3 = (-)^{j-1} i \pt [\pt \Phi \pt dz\wedge d\bar{z} + \Phi^{-1}(du + A \hp dz) \wedge (d\bar{u} + \bar{A} \hp d\bar{z})] \rlap{.} \label{dimH=1-GLT-sympl-forms}
\end{gather}
For $j \geq 2$, by expressing the inverse matrices which occur in the formulas for $\Phi$ and $A$ in terms of adjugate matrices and thus of minors and then judiciously repositioning in the second case the result so as to be able to use Laplace's cofactor expansion formula, we can show that these can be expressed entirely in terms of determinants of Hankel matrices.\footnote{\,Hankel matrices are matrices whose entries along the parallels to the secondary diagonal are equal for each parallel in part.} Specifically, we get 
{\allowdisplaybreaks
\begin{equation}
\begin{aligned}
\Phi & = \phantom{- i \hp} \frac{\det(\bar{h}_{m+m'})_{1-j \leq m,m' \leq j-1} \hfill}{\det(\bar{h}_{m+m'})_{3/2-j \leq m,m' \leq j-3/2}} \\
A & = - i \hp \frac{\det(\bar{h}_{m+m'})_{3/2-j \leq m,m' \leq j-1/2}}{\det(\bar{h}_{m+m'})_{3/2-j \leq m,m' \leq j-3/2}} \rlap{.}
\end{aligned}
\end{equation}
}%
Further manipulations of the determinants involving the multiplication of some rows or columns by a sign yield finally the expressions
{\allowdisplaybreaks
\begin{equation} \label{dimH=1-GLT-HA}
\begin{aligned}
\Phi & = \phantom{i \hp} \frac{\det(h_{m+m'})_{1-j \leq m,m' \leq j-1} \hfill}{\det(h_{m+m'})_{3/2-j \leq m,m' \leq j-3/2}} \\
A & = i \hp \frac{\det(h_{m+m'})_{1/2-j \leq m,m' \leq j-3/2}}{\det(h_{m+m'})_{3/2-j \leq m,m' \leq j-3/2}} \rlap{.}
\end{aligned}
\end{equation}
}%
These formulas generalize to arbitrary values of $j \geq 2$ the $j=2$ case formulas found in \cite{MR2177322}.

\subsection{}

The arbitrariness implicit in the singling out of a particular complex structure is in some sense the most unsatisfactory feature of the formulas in Proposition \ref{GLT-formulas}. We know that in the $j=1$ case there exist expressions for the hyperk\"ahler metric and symplectic forms in which symmetries rather than the complex structure play the more prominent role, and these are of course the Gibbons-Hawking formulas. To see how the Gibbons-Hawking Ansatz arises in the Legendre transform approach we need to perform an additional change of variables. The basic idea is to try to cast as coordinates the coefficients $x^{\I}_m$\,---\,or, rather, their real components, which we define by $x^{\I}_{\pm} = \pm \frac{1}{2} (x^{\I}_1 \pm i x^{\I}_2)$ and $x^{\I}_0 = x^{\I}_3$. These can be viewed as forming the Euclidean components of some $\mathbb{R}^3$-vectors $\vec{r}^{\,\I}$. This, however, does not give us enough coordinates, since we have only three of them for each value of the index $I$. To complete the number needed one notices that the coefficients $x^{\I}_m$\,---\,implicitly expressed in terms of $z^{\I},u_{\I}$ and their complex conjugates by means of the Legendre transform relations \eqref{LT-eqs} (recall also that for $j=1$, $x^{\I}_{+}=z^{\I}$)\,---\,are invariant with respect to shifts along the real directions of the complex $u_{\I}$-variables. These directions can then be parametrized to furnish a fourth real coordinate for each value of $I$, which we denote by $\psi_{\I}$. We have thus
\begin{equation}
u_{\I} = \psi_{\I} + \frac{i}{2}L_{x^I} \rlap{.}
\end{equation}
In terms of the new variables the formulas of Proposition \ref{GLT-formulas} take then the familiar form
{\allowdisplaybreaks
\begin{equation}
\begin{aligned}
g & = \frac{1}{2} U_{\I\J}\pt d\vec{r}^{\,\I} \!\cdot d\vec{r}^{\,\J} + \frac{1}{2} U^{\I\J}(d\psi_{\I} + A_{\I}) (d\psi_{\J} + A_{\J}) \\
\vec{\omega} & = - \frac{1}{2}\pt U_{\I\J}\pt d\vec{r}^{\,\I} \!\wedge d\vec{r}^{\,\J} - d\vec{r}^{\,\I} \!\wedge (d\psi_{\I} + A_{\I}) 
\end{aligned}
\end{equation}
}%
with matrix inverses indicated by upper indices and the wedge product of two $\mathbb{R}^3$-vector-valued 1-forms defined as one would expect from vector calculus, where we identify
\begin{equation} \label{GH-LT-eqs}
U_{\I\J} = - \frac{1}{2} L_{x^{\I}_0x^{\J}_0}
\quad\text{and}\quad
A_{\I} = \Im (L_{x^{I}_0x^{J}_{+}}dx^{\J}_{+}) \rlap{.}
\end{equation}
This metric has an obvious tri-Hamiltonian $\mathbb{R}^n$-symmetry generated by shifts in the $\psi_{\I}$ directions, with corresponding moment map images $\vec{r}^{\,\I}$. Finally, the Bogomolny field equations satisfied by the Higgs fields $U_{\I\J}$ and connection 1-forms $A_{\I}$ can be understood as consequences of the differential equations satisfied by $L$. 

This concludes our review of Lindstr\"om and Ro\v{c}ek's generalized Legendre transform construction. 

\bigskip

\textit{N.B. In the remainder of the paper we will be concerned exclusively with four-dimensional hyperk\"ahler metrics and consequently, since in this case they take only one value, we will automatically drop the indices $I,J,\dots$~from the formulas and notations of this section whenever we resort to them.}

\section{Majorana polynomials} \label{sec:Majorana-pols}

\subsection{}

A straightforward attempt to pursue in the cases with $j \geq 2$ the same strategy of casting the coefficients $x_m$ as coordinates which has worked in the $j=1$ case runs instantly into fatal difficulties because the number of these, $2j+1$, is higher than the four which are needed for dimensional reasons. The generalized Legendre transform equations \eqref{GLT-eqs} make it otherwise very clear that some of them need to be eliminated, so any na\"ive generalization in this direction is condemned from the start. 

So then perhaps we could try to look at the $j=1$ case through another pair of glasses in a way which is more amenable to generalization. This turns out to be indeed possible provided we put on our spherical coordinates glasses. As we will see in the next section, the Gibbons-Hawking Ansatz, in its full generality and not just when a rotational symmetry is present, admits a natural description in terms of spherical coordinates. This is, more than anything, a reflection of the quaternionic structure of the hyperk\"ahler space, and so it is actually more profitable and indeed natural to use quaternions instead of spherical angles. 

Another clue in this direction comes from Quantum Mechanics. Polynomial functions of degree $2j$ formally identical to the local sections
\begin{equation} \label{Majorana-pol}
x_N(\zeta) = \sum_{m=-j}^{j} x_m \hp \zeta^{j-m}
\end{equation}
have been used as far back as 1932 by E. Majorana to describe wave functions of quantum-mechanical particles with spin $j$ \cite{Majorana1932}. In this description generalizing the Bloch sphere representation of particles with spin $1/2$, states with spin $j$ are represented by $2j$ points on the Riemann sphere, usually referred to in this context as the Majorana sphere, corresponding to the roots of this polynomial. In our case, the reality condition
\begin{equation} \label{x-alt-real}
\bar{x}_m = (-)^mx_{-m}
\end{equation}
satisfied by the polynomial coefficients implies that the $2j$ points come in anti\-po\-dally-opposite pairs. The twistor-theoretic potential of Majorana's ideas was recognized by R. Penrose,  see \textit{e.g.} \cite{MR1048125} and \cite{MR1865778} (esp.~Appendix C). 

\subsection{}

Mathematically, this points at the spin-$j$ representations of the group $SU(2)$. To cast the connection in more precise terms it helps to recall first a few basic facts about these. In line with our considerations above it is useful to view in this context the Lie group $SU(2)$ in its concrete manifestation as the group of unitary quaternions. An explicit unitary $2j+1$-dimensional matrix realization of this is given by the Wigner D-matrices. To any non-zero quaternion $q = q_0 + q_1 \pt \pmb{i} + q_2 \pt \pmb{j} + q_3 \pt \pmb{k} \in \mathbb{H}^{\times}$ one associates a matrix (see \textit{e.g.} \cite{hanson2006})
\begin{align}
& D^j_{mm'}(q) = \frac{\sqrt{(j+m)!(j-m)!(j+m')!(j-m')!}}{(|v|^2+|w|^2)^j} \\
& \qquad \times i^{-2j} \sum_{s} \frac{w^s}{s!} \frac{\bar{w}^{m-m'+s}}{(m-m'+s)!} \frac{v^{j+m'-s}}{(j+m'-s)!} \frac{(-\bar{v})^{j-m-s}}{(j-m-s)!} \nonumber
\end{align}
where $m,m' = -j,\dots,j$ and by definition $w = q_3 + i\hp q_0$ and $v = q_1+i\hp q_2$. In general $j$ can take positive integer as well as half-integer values, although here we will be concerned only with the former case. The essential property of these matrices is that they provide a representation of the multiplicative group of non-zero quaternions $\mathbb{H}^{\times}$, that is, $D^j(q)D^j(q') = D^j(q\hp q')$ for all $q,q' \in \mathbb{H}^{\times}$. Given their explicit form one can easily verify that they satisfy also the following properties:
\begin{alignat}{4} \label{Wigner-D-propr}
1. \ & \text{Scale invariance:} \quad	&& D^j_{mm'}(q/|q|) = D^j_{mm'}(q) \nonumber \\[3pt]
2. \ & \text{Reflection:}             	&& D^j_{mm'}(-q) = (-)^{2j} D^j_{mm'}(q) \\
3. \ & \text{Time-reversal:} 		&& \overline{D^j_{mm'}(q)} = (-)^{m-m'} D^j_{-m,-m'}(q) \rlap{.} \nonumber
\end{alignat} 
The first of these shows that the representation descends in fact to the subgroup of \textit{unit} quaternions, which is isomorphic to $SU(2)$. The second property shows that for integer values of $j$ the representation descends further to the group $SO(3)$ for which $SU(2)$ is a double cover. The name of the third property derives from a certain role it plays in Quantum Mechanics.

As anyone familiar with them knows, Wigner D-matrices are usually given in the literature in terms of Euler angles rather than quaternions. The link between the two parametrizations is provided by the relations
\begin{equation}
w = |q| \sin \frac{\theta}{2} \pt e^{\frac{i}{2}(\phi - \psi)}
\qquad\qquad
v = |q| \cos \frac{\theta}{2} \pt e^{\frac{i}{2}(\phi + \psi)} \rlap{.}
\end{equation}
From a computational point of view it is important to note that the quaternionic expressions are vastly superior because they involve only homogeneous polynomials in the quaternion's components. 

\subsection{}

Rotations of the Majorana sphere translate into unitary linear transformations of the $2j+1$ coefficients $x_m$ in accordance with the spin-$j$ representation of the rotation group. More precisely, the \textit{normalized} coefficients transform in this way, where we introduce a numerical spherical tensor normalization factor for each coefficient by $x^j_m \equiv c^j_{m} x^{\phantom{j}}_m$, with
\begin{equation} \label{sf-norm}
c^j_m = \sqrt{\frac{(j+m)!(j-m)!}{(2j)!}} \rlap{.}
\end{equation}
Then such transformations are of the form 
\begin{equation} \label{spinor-rot}
x^j_m = \sum_{m'=-j}^j D^j_{mm'}(q) \pt \mathring{x}^j_{m'}
\end{equation}
with the quaternion $q$ parametrizing the rotation and $\mathring{x}^j_m$ indicating the (normalized) initial coefficients, before the rotation. The third property \eqref{Wigner-D-propr} of the Wigner D-matrices guarantees that the reality property \eqref{x-alt-real} of the coefficients is preserved. For $j=1$ this is the usual transformation law of a 3-\textit{vector} expressed in a complex spherical basis, while for $j \geq 2$ it is the transformation law of a \textit{spinor}. For a  generic $j$, we will use the term ``spinor" in an extended sense to refer to either case.

But rather than pursue an active point of view of these transformations we prefer to assume a passive one and use them to introduce a spherical or, which is the same thing in our approach, a quaternionic parametrization of the spinor. This can be done by fixing a \textit{reference spinor} $\mathring{x}^j_m$ and then parametrizing every spinor in the orbit of the unitary transformation \eqref{spinor-rot} by the parameters of the specific Wigner matrix used to reach it from this. 

The choice and parametrization of the reference spinor itself is in general a non-trivial matter. First, note that a counting of the degrees of freedom implies that the reference spinor needs to depend on $2j-2$ parameters: the $2j+1$ parameters of the $x^j_m$'s minus the three parameters of unitary rotations (when counting rotational parameters, Euler angles rather than quaternions are a better choice). The case $j=1$ is an exception to this rule because to spherically-parametrize a 3-vector one needs only two angles, not three, which leaves the reference vector depending on one parameter, not zero. This is of course the length of the vector, which is an invariant at rotations. In fact, the same principle applies for higher values of $j$, where the reference spinor should be parametrized entirely by $2j-2$ independent rotation invariants. Apart from this constraint, the choice of reference spinor is then arbitrary. And we mean this in the same way as, say, choosing the position of the north pole when spherically-parametrizing the Earth is arbitrary. Although, as in this example, sometimes extraneous criteria may guide one's choices as well. 

So then when choosing a reference spinor one is faced with the problem of constructing spherical invariants. One way to do this is by means of spherical tensor couplings. The two lowest-order invariants that can be obtained in this way are
\begin{equation} \label{rot-invs}
\begin{gathered}
r_2(x) = \frac{1}{2} \hspace{-1pt} \sum_{m=-j}^j |x^j_{m}|^2 \\[-2pt]
r_3(x) = \frac{N_j}{3!} 
\hspace{0pt} 
\sum_{\shortstack{$\scriptstyle m_1,m_2,m_3=-j$ \\ $\scriptstyle m_1+m_2+m_3=0$}}^j
\hspace{-2pt} 
\left(
\begin{array}{ccc}
\! j & \!\!\! j & \!\!\! j \! \\
\! m_1 & \!\!\! m_2 & \!\!\! m_3 \!
\end{array}
\right)
x^j_{m_1} x^j_{m_2} x^j_{m_3}
\end{gathered}
\end{equation}
where the array in the second formula indicates a Wigner $3j$-symbol and $N_j$ is an arbitrary numerical normalization constant. Their invariance is manifest by construction and one can indeed check that $r_2(x) = r_2(\mathring{x})$ and $r_3(x) = r_3(\mathring{x})$. Notice that the second one vanishes trivially for $j=1$. On the other hand, for higher values of $j$ one can construct further non-trivial higher-order invariants by similar means. These may or may not be independent, a thing one needs to check. 

After constructing in this manner the required number of independent invariants we are still left with the obvious problem posed by the fact that their scaling properties are not what we need them to be (they have to scale the same as the $x^j_m$'s). To obtain a good set of invariant parameters for the reference spinor we need to further solve an algebraic equation with these invariants as coefficients. Rather than continue the discussion in the generic case we refer the reader to the subsequent two sections where we work this out in detail for the cases with $j=1$ and $j=2$. 

\subsection{}

Let us assume now that we have managed to choose a good reference spinor\,---\,where, we stress once again, the equation \eqref{spinor-rot} has to be viewed simply as a change of variables, a reparametrization of the spinor components. As usual in such cases, in applications we will need the Jacobian matrix of the transformation, which is to say, the induced differential map. Differentiating the equation \eqref{spinor-rot} we obtain
\begin{equation} \label{diff-map}
dx^j_m = \sum_{m'=-j}^j D^j_{mm'}(q) \pt \chi^j_{m'}
\end{equation}
with
\begin{equation} \label{chi-basis}
\chi^j_m = d\mathring{x}^j_m \, + \, 2\hp i \sqrt{2j(j+1)} 
\hspace{-5pt} 
\sum_{\shortstack{$\scriptstyle m_1,m_2$ \\ $\scriptstyle m_1+m_2=m$}} 
\hspace{-5pt} 
\left(
\begin{array}{c|cc}
\! j & 1 & \!\!\! j \! \\
\! m & m_1 & \!\!\! m_2 \!
\end{array}
\right)
\sigma^1_{m_1} \mathring{x}^j_{m_2} \rlap{.}
\end{equation}
These formulas allow us to read off immediately the differential map for any value of $j$ and any choice of reference spinor $\mathring{x}^j_m$. The array symbol denotes a Clebsch-Gordan coefficient and the $\sigma^1_n$ with $n=-1,0,+1$ are normalized left-invariant 1-forms for the $SU(2)$ group. More precisely, let $q^{-1} dq = \sigma_0 + \sigma_1 \pt \pmb{i} + \sigma_2 \pt \pmb{j} + \sigma_3 \pt \pmb{k}$. The quaternionic algebra structure implies that the real-valued component 1-forms satisfy the Cartan-Maurer relations $d\sigma_0 = 0$, $d\sigma_1 = - 2 \pt \sigma_2 \wedge \sigma_3$, a.s.o.~A spherical basis is obtained by the complex-linear transformation $\sigma_{\pm} = \mp \frac{1}{2} (\sigma_1 \mp i \hp \sigma_2)$, $\sigma_0 = \sigma_3$, and we then normalize by $\sigma^1_n = c^1_n\sigma^{\phantom{1}}_n$, where the normalization coefficients are given as before by the combinatoric formula \eqref{sf-norm}. The reality property $\bar{\sigma}_{n} = (-)^n \sigma_{-n}$ together with that of the reference spinors implies that the frame elements $\chi^j_m$ satisfy as well an alternating reality property. 

\subsection{}

We end this section by describing two alternative representations of Majorana polynomials which will prove better suited than the one we have discussed so far for certain arguments that we will need to make. As we have mentioned already, the reality property \eqref{x-alt-real} satisfied by the coefficients of Majorana polynomials implies that the polynomial roots come in antipodally-conjugated pairs. By way of the fundamental theorem of Algebra we can show that such polynomials can always be cast into the following canonical factorized form
\begin{equation} \label{Majorana-roots}
x_N(\zeta) = \rho \prod_{i=1}^j \frac{(\zeta-a_i)(1+\bar{a}_i\zeta)}{(1+|a_i|^2)}
\end{equation}
with the overall scaling factor $\rho$ real. This representation is formally very closely related to the so-called \textit{spin coherent state} representation of wave functions in Quantum Mechanics. Notice that the way we assign the labels $a_i$ and $a^c_i$ to a pair of antipodally-conjugated roots is really arbitrary\,---\,and, moreover, if we switch the labeling, then the corresponding factor in the product changes sign. Because of this fact we can always label the roots in such a way as to have $\rho > 0$. 

If we write $a_i = v_i/\bar{w}_i$, for complex numbers $v_i,w_i$ defined up to a real scale, which we then fix by requiring that $\rho = \prod_{i=1}^j(|v_i|^2 + |w_i|^2)$, we obtain the equivalent representation
\begin{equation} \label{wv-Majorana}
x_N(\zeta) = \prod_{i=1}^j (w_i + \zeta\bar{v}_i)(v_i - \zeta\bar{w}_i) \rlap{.}
\end{equation}
The pairs $w_i,v_i$ can be regarded in some sense as the complex components of $j$ quaternions.

\section{$\mathcal{O}(2)$ constructions in a quaternionic frame}

In the $j=1$ case the presence of the abelian tri-Hamiltonian toric symmetry provides us with a very special symmetry-adapted coordinate system, the Gibbons-Hawking coordinates. For higher values of $j$ the symmetries take on a hidden character and, as we have noted earlier, simplistic attempts at generalizing the Gibbons-Hawking frame are readily proven unfeasible. The idea then is to try to look at the $j=1$ case  in a way which may not be as immediately gratifying as the Gibbons-Hawking approach, but which on the other hand has a better potential for generalization. In the previous section we have already started making the case that such a way is offered by spherical or, equivalently, quaternionic coordinate frames. In this section we revisit the $j=1$ case from this point of view, as a warm-up exercise for the more challenging $j=2$ case which we will discuss in the next section. 

From a twistorial perspective the case with $j=1$ is characterized by the factorization of the twistor space holomorphic projection onto the twistor sphere through an intermediate $\mathcal{O}(2)$ bundle (as stated earlier, all of our considerations here will be restricted to four-dimensional hyperk\"ahler spaces). This $\mathcal{O}(2)$ bundle has a globally-defined section characterized in a certain trivialization by three moduli denoted by $x_m$, with $m=-1,0,+1$. Alternatively, one can view these as the complex spherical basis components of a real Euclidean 3-vector $\vec{r}$, and indeed, rotations of the twistor sphere induce on them linear transformations in accordance with the spin-1 or vector representation of the $SO(3)$ group. Geometrically, the vector $\vec{r}\hp$~is interpreted as the image of the hyperk\"ahler moment map for the toric action with respect to the standard frame in the space of hyperk\"ahler symplectic forms, and in the Gibbons-Hawking approach is used together with a toric orbit coordinate to coordinatize the space.

For $j=1$, there is only one non-trivial $SO(3)$ invariant that one can construct, namely the second-order invariant in \eqref{rot-invs}, with higher-order invariants being either trivially zero or dependent on this. In terms of the unnormalized components it can be written as
\begin{equation}
r_2(x) = \frac{1}{2} \pt \big[ \pt |x_{+1}|^2 + \frac{1}{2}(x_0)^2 + |x_{-1}|^2 \hp \big] = -
\left|
\begin{array}{cc}
\hspace{-2pt} \displaystyle \frac{x_{+1}}{1} & \displaystyle \frac{x_{0}}{2}  \hspace{-2.5pt} \\[7pt]
\hspace{-2pt} \displaystyle \frac{x_{0}}{2} & \displaystyle \frac{x_{-1}}{1} \hspace{-2.5pt}
\end{array}
\right|
\end{equation}
where the second determinant expression is a wink in the direction of the next section. This invariant is of course, up to a numerical factor, nothing but the Euclidean length of the 3-vector $\vec{r}$ (squared). As reference vector we choose
\begin{equation}
(\hp \mathring{x}_m \hp )_{m = -1, 0, +1} =
\begin{pmatrix}
0 \\[0pt]
\, \rho \,  \\[2pt]
0
\end{pmatrix}
\end{equation}
and so we have $r_2(x) = r_2(\mathring{x}) = \rho^2/4$. We take then $\rho$ to be by definition the positive solution to this quadratic equation. One can easily see that this is in fact the same $\rho$ as the one which occurs in the $j=1$ version of formula \eqref{Majorana-roots}.

With this choice, the passage to spherical/quaternionic variables induces a differential map described by formula \eqref{diff-map}, with $j=1$ and
\begin{equation} \label{chi-basis-j=1}
\begin{aligned}
\chi_{+1} & = -i \hp \rho \pt (\sigma_1-i\hp\sigma_2) \\[1pt]
\chi_{\hp 0 \phantom{+}} & = d\rho \\[0pt]
\chi_{-1} & = -i \hp \rho \pt (\sigma_1+i\hp\sigma_2) \rlap{,}
\end{aligned}
\end{equation}
where the spherical normalization factors have been stripped off.

Denoting $\bar{h}_{m} = L_{x_{k}x_{m-k}}$ in accordance with \eqref{h-bar-def} and normalizing these by means of the \textit{dual} normalization convention $h^1_m \equiv (c^1_{m})^{-1} h^{\phantom{1}}_m$, we notice that in quaternionic dimension one the Gibbons-Hawking connection 1-form \eqref{GH-LT-eqs} can be rewritten as
{\allowdisplaybreaks
\begin{align}
A & = - \frac{i}{\sqrt{2}} 
\sum^1_{\shortstack{$\scriptstyle m_1,m_2 = -1$ \\ $\scriptstyle m_1+m_2=0$}} 
\hspace{-3pt} 
\left(
\begin{array}{c|cc}
\! 1 & 1 & \!\!\! 1 \! \\
\! 0 & m_1 & \!\!\! m_2 \!
\end{array}
\right)
h^1_{m_1} dx^1_{m_2} \\[-1pt]
& = - \frac{i}{\sqrt{2}}
\sum^1_{\shortstack{$\scriptstyle m,m_1,m_2 = -1$ \\ $\scriptstyle m_1+m_2=m$}} 
\hspace{-5pt}
D^1_{0m}(q) 
\left(
\begin{array}{c|cc}
\! 1 & 1 & \!\!\! 1 \! \\
\! m & m_1 & \!\!\! m_2 \!
\end{array}
\right)
\mathring{h}^1_{m_1} \chi^1_{m_2} \rlap{.} \nonumber
\end{align}
}%
The second line follows from the equivariance properties of the Clebsch-Gordan coefficients, and by definition we set
\begin{equation}
\mathring{h}^1_m = \sum_{m'=-1}^1 D^1_{mm'}(q^{-1}) \hp h^1_{m'} \rlap{.}
\end{equation}

Finally, let us observe that if we introduce through a further linear transformation the complex-valued 1-forms
\begin{equation} \label{theta-basis-j=1}
\begin{aligned}
\theta_{+1} & = h_{0} \hp \chi_{+1} \\[-2pt]
\theta_{\hp 0\phantom{+}} & = \frac{1}{2} \hp h_{0} \hp \chi_{\hp 0} + i \hp (d\psi + A) 
\end{aligned}
\end{equation}
then in terms of them the Gibbons-Hawking metric and symplectic forms assume the following manifestly quaternionic expressions
{\allowdisplaybreaks
\begin{align} \label{HKstr-j=1}
g = \lambda \pt (\pt |\theta_{\hp 0}|^2 + |\theta_{+1}|^2) \qquad\text{and}\qquad & \omega^1_m =  \sum_{m'=-1}^1 D^1_{mm'}(q) \pt \mathring{\omega}^1_{m'} \\
& \text{with} \nonumber \\[2pt]
& \mathring{\omega}_{+1} =  \lambda \hp i \, \theta_{+1} \wedge \bar{\theta}_{\hp 0} \nonumber \\[0pt]
& \mathring{\omega}_{\hp 0} \hspace{6.05pt} = \lambda \hp i \pt (\hp \theta_{\hp 0} \wedge \bar{\theta}_{\hp 0} - \theta_{+1} \wedge \bar{\theta}_{+1}) \nonumber \\[0pt]
& \mathring{\omega}_{-1} = \lambda \hp i \, \bar{\theta}_{+1} \wedge \theta_{\hp 0} \rlap{,} \nonumber
\end{align}
}%
where we spherically-normalize the latter by $\omega^1_m = c^1_m \omega^{\phantom{1}}_m$ and we have
\begin{equation}
\lambda = - \frac{1}{h\rlap{$_{0}$}} \rlap{\,.}
\end{equation}

Thus, to sum up, we have shown that the $j=1$ hyperk\"ahler metric and symplectic forms
\begin{itemize}
\setlength\itemsep{1pt}
\item[1.] are entirely determined by the three functions $h_m$ with $m=-1,0,1$;
\item[2.] can be re-cast in a manifestly quaternionic canonical form by means of two successive transformations, the coordinate transformation \eqref{diff-map}--\eqref{chi-basis-j=1} and the linear change of frame \eqref{theta-basis-j=1}. 
\end{itemize}

\section{$\mathcal{O}(4)$ constructions in a quaternionic frame} \label{sec:O4-constr}

In the $j=2$ case we have no correspondent of the Gibbons-Hawking Ansatz to start with\,---\,finding one such correspondent in some sense is in fact precisely what we are after\,---\,but we do possess nonetheless expressions for the hyperk\"ahler metric and symplectic forms, namely the generalized Legendre transform ones given in Proposition \ref{GLT-formulas} and the discussion following it. These are written in a holomorphic coordinate frame associated to a particular hyperk\"ahler complex structure, and our aim here is to reformulate them in a quaternionic coordinate frame. 

\subsection{}

The holomorphic twistor space projection to the twistor $\mathbb{CP}^1 \cong S^2$ factorizes now through an $\mathcal{O}(4)$ bundle possessing a global section characterized in our standard trivialization by the five moduli $x_m$ with $m=-2,\dots,2$ satisfying an alternating reality condition. At rotations of the twistor $S^2$ these transform according to the unitary spin-2 representation of $SO(3)$ and thus can be thought of as the components of a spinor. From them one can construct this time not one but two (and only two) independent non-trivial real-valued $SO(3)$ invariants, which we may write as follows:
{\allowdisplaybreaks
\begin{equation} \label{r2r3}
\begin{gathered}
r_2(x) = \frac{1}{2} \pt \big[ \pt |x_{+2}|^2 + \frac{1}{4}|x_{+1}|^2 + \frac{1}{6}(x_0)^2 + \frac{1}{4}|x_{-1}|^2 + |x_{-2}|^2 \big] \\[3pt]
r_3(x) =
\left|
\begin{array}{ccc}
\hspace{-2pt} \displaystyle \frac{x_{+2}}{1} & \displaystyle \frac{x_{+1}}{4} & \displaystyle \frac{x_0}{6} \hspace{-2.5pt} \\[7pt]
\hspace{-2pt} \displaystyle \frac{x_{+1}}{4} & \displaystyle \frac{x_0}{6} & \displaystyle \frac{x_{-1}}{4} \hspace{-2.5pt} \\[7pt]
\hspace{-2pt} \displaystyle \frac{x_0}{6} & \displaystyle \frac{x_{-1}}{4} & \displaystyle \frac{x_{-2}}{1} \hspace{-2.5pt}
\end{array}
\right| \rlap{.}
\end{gathered}
\end{equation}
}%
Indeed, inserting spherical normalizations, one can easily verify that these two expressions coincide precisely with the two manifestly invariant expressions \eqref{rot-invs} corresponding to $j=2$ provided that in the second case we choose the normalization constant to be $N_2 = \sqrt{35/12}$. For reasons which will soon become clear we actually prefer to work with the two scaled versions
\begin{equation} \label{Weierstrass-coeffs}
\begin{aligned}
g_2 & = 4 \pt r_2(x) \\[1pt]
g_3 & = 16 \pt r_3(x) \rlap{.}
\end{aligned}
\end{equation}

Before we address the issue of choosing a reference spinor we wish to open a parenthesis to make an observation whose usefulness will become apparent only later on, concerning the spherical coupling process through which the invariants have been constructed. Note that similarly to how we construct invariants we can employ spherical coupling to construct spherical tensors, that is, quantities which transform covariantly with respect to $SO(3)$ rotations in accordance with one of its finite-dimensional representations. However, in some interesting cases we can achieve the same result by differentiating the invariants. Thus, consider the quantities defined by
{\allowdisplaybreaks
\begin{equation} \label{dgdx}
\begin{aligned}
d_{2,m} & = (-)^m \frac{\partial g_2}{\partial x_{-\rlap{$\scriptstyle m$}}} \\[-1pt]
d_{3,m} & = (-)^m \frac{\partial g_3}{\partial x_{-\rlap{$\scriptstyle m$}}} \rlap{\ .}
\end{aligned}
\end{equation}
}%
These satisfy alternating reality properties and, for $i=2,3$, introducing dual normalizations by $d^{\pt 2}_{i,m} = (c^2_m)^{-1}d^{\phantom{2}}_{i,m}$, one can check that at transformations \eqref{spinor-rot} with $j=2$ they transform covariantly according to the rule
\begin{equation}
d^{\pt 2}_{i,m} = \sum_{m'=-2}^2 D^2_{mm'}(q) \pt \mathring{d}^{\pt 2}_{i,m'} \rlap{.}
\end{equation}

Returning now to our discussion, to choose a reference spinor we make the Ansatz 
\begin{equation} \label{ref-sp-j=2}
(\hp \mathring{x}_m \hp )_{m=-2,\dots,+2} = \frac{1}{4}
\begin{pmatrix}
e_1-e_3 \\
0 \\
\! 6 \hp (e_1+e_3) \! \\
0 \\
e_1-e_3
\end{pmatrix} \rlap{.}
\end{equation}
Introducing a further dependent parameter $e_2$ such that $e_1+e_2+e_3 = 0$, by invariance we have then $g_2 = 4 \pt r_2(\mathring{x}) = - (e_1e_2+e_2e_3+e_3e_1)$ and $g_3 = 16 \pt r_3(\mathring{x}) = e_1e_2e_3$. In other words, what we have shown is that \eqref{ref-sp-j=2} represents a good choice of reference spinor provided that $e_1,e_2,e_3$ are the three roots of the cubic Weierstrass polynomial formed with $g_2$ and $g_3$, that is
\begin{equation}
X^3 - g_2X - g_3 = (X-e_1)(X-e_2)(X-e_3) \rlap{.}
\end{equation}

The Weierstrass discriminant $\Delta = 4 g_2^3 - 27g_3^2 = (e_1-e_2)^2(e_1-e_3)^2(e_2-e_3)^2$ is, very importantly, not only real but also positively defined. This is best seen if we resort to the representation \eqref{wv-Majorana} of the Majorana polynomial corresponding to $j=2$. Comparison with the defining representation \eqref{Majorana-pol} yields by way of Vieta's relations expressions for the coefficients $x_m$ in terms of the variables $w_i,v_i$ and their complex conjugates. Substituting these into the formulas for $g_2$, $g_3$ and then the resulting expressions into that of the discriminant yields eventually $\Delta = [\pt (|v_1|^2 + |w_1|^2) (|v_2|^2 + |w_2|^2) |v_1\bar{v}_2 + w_2\bar{w}_1|^2 |v_1\bar{w}_2 - v_2\bar{w}_1|^2 \pt ]^2$, which is manifestly positive. 

Positivity of the Weierstrass discriminant means that the Weierstrass roots are real. As such, they can be ordered along the real axis, and we adopt here the usual ordering convention $e_1 > e_2 > e_3$. The method above can give us in fact for each factor of $\Delta$ the expressions
{\allowdisplaybreaks
\begin{equation}
\begin{aligned}
e_1 - e_2 & = |v_1\bar{v}_2 + w_2\bar{w}_1|^2 \\
e_1 - e_3 & = (|v_1|^2 + |w_1|^2) (|v_2|^2 + |w_2|^2) \\
e_2 - e_3 & = |v_1\bar{w}_2 - v_2\bar{w}_1|^2 \rlap{.}
\end{aligned}
\end{equation}
}%
In particular, this shows that we can identify $e_1 - e_3 = \rho$, as defined in equation \eqref{Majorana-roots} (with $j=2$). 

The upshot of these considerations is the following: the $j=2$ version of equation \eqref{spinor-rot} together with this choice of reference spinor provide a parametrization of the five-component spinor $x_m$ in terms of three spherical parameters, the Euler angles associated with the quaternion $q$, plus two spherical invariants, $e_1$ and $e_3$ (or, equivalently, and this will be in fact our preferred choice, $\rho$ and $e_2$). This is in a very clear sense a spinorial generalization of the usual spherical parametrization of $3$-vectors. 

The differential map corresponding to the reparametrization is described by the equation \eqref{diff-map}, with equation \eqref{chi-basis} giving
\begin{equation} \label{chi-basis-j=2}
\begin{aligned}
\chi_{+2} & = \frac{1}{4} \hp d\rho - i \hp \rho \pt \sigma_3 \\[2pt]
\chi_{+1} & = - i \hp \rho \pt (\sigma_1 +i \hp \sigma_2) + 3\hp i \hp e_2 (\sigma_1-i \hp \sigma_2) \\[-2pt]
\chi_{\hp 0\phantom{+}} & = - \frac{3}{2} \hp d e_2 \rlap{\pt .}
\end{aligned}
\end{equation}
The remaining elements are obtained by (alternating) conjugation. 

\subsection{}

Our discussion so far has centered on a change of variables and has not touched yet on the generalized Legendre transform equations. Let us turn now to these and examine them closer. For $j=2$ and in quaternionic dimension one they read: 
\begin{equation} \label{GLT-eqs-j=2}
L_{x_m} = 
\begin{cases}
- i u & \text{for $m=+1$} \\
\hspace{2.2ex} 0 & \text{for $m=0$}  \\
-i \bar{u} & \text{for $m=-1$} \rlap{.}
\end{cases}
\end{equation}
It is understood that $L$ is a (possibly transcendental) function of the five variables $x_m$. The main difference with respect to the $j=1$ case is the presence of the differential constraint $L_{x_0} = 0$. This constraint is the source of all difficulties and the crux around which the whole matter revolves. In principle it allows us to express implicitly, say, $x_0$ in terms of the other four variables. Then the remaining two equations \eqref{GLT-eqs-j=2} corresponding to $m = \pm 1$ together with the identification $x_{+2} = z$  let us further exchange, again implicitly, these four variables for the four holomorphic coordinates $z,u,\bar{z},\bar{u}$\,---\,in the differential frame associated to which we have expressions for the hyperk\"aler metric and symplectic forms. 

Here, however, we propose a different approach. Our key idea is to employ instead of the five variables $x_m$ the equivalent set of spherical parameters $\phi$, $\theta$, $\psi$, $\rho$ and $e_2$, and then use the constraint $L_{x_0} = 0$ to solve implicitly for $e_2$ in terms of $\phi$, $\theta$, $\psi$ and $\rho$.\footnote{\hp Note once again that here and throughout the paper we \textit{say} Euler angles but we \textit{mean} quaternions. The former are better for variable counting purposes while the latter offer calculational ease and conceptual clarity.} The latter reduced set gives us, very importantly, a system of coordinates on the hyperk\"ahler space. So, rather than pass to the holomorphic coordinate frame, we choose to transition the other way around and express the hyperk\"aler metric and symplectic forms in the spherical (or, better, \textit{quaternionic}) frame given by $\sigma_1, \sigma_2, \sigma_3$ and $d\rho$. 

Concretely, we carry this out as follows: differentiating the generalized Legendre relations and denoting as in \eqref{h-bar-def} the double derivatives $L_{x_{m}x_{\smash{m'}}} = \bar{h}_{m+m'}$, we obtain
\begin{equation} \label{diff-GLT-eqs}
\sum_{m'=-2}^2 \bar{h}_{m+m'}dx_{m'} = 
\begin{cases}
-i \hp du & \text{for $m=+1$} \\
\hspace{2.7ex} 0 & \text{for $m=0$} \\
-i \hp d\bar{u} & \text{for $m=-1$} \rlap{.}
\end{cases}
\end{equation}
We then proceed in succession to substitute into the $m=0$ component of this equation the corresponding differential map relations \eqref{diff-map}, solve the result for $\chi_0$ in terms of $\chi_{+2}, \chi_{+1}, \chi_{-1}, \chi_{-2}$, and then substitute this back into the relations \eqref{diff-map}. Note that since the equations \eqref{chi-basis-j=2} can be neatly separated into $\chi_0 = - \frac{3}{2} \hp d e_2$ on one hand and
\begin{equation} \label{quat-chi-frame}
\begin{pmatrix}
\chi_{+2} \\[1pt]
\chi_{+1} \\[1pt]
\chi_{-1}  \\[1pt]
\chi_{-2}
\end{pmatrix}
=
\left(
\hspace{-5pt}
\begin{array}{c@{\ }c@{\ \, }c@{\ \ }c@{}}
\, 1/4		& 0 					& 0 						& \ \llap{$-$} i \rho \ 	\\
\, 0 		& 3\hp i\hp e_2 - i\hp \rho 	& \phantom{+} 3\hp e_2 + \rho 	& 0 \				\\
\, 0 		& 3\hp i\hp e_2 - i\hp \rho 	& -3\hp e_2 - \rho			& 0 \				\\
\, 1/4		& 0 					& 0 						& i \rho \
\end{array}
\hspace{-1pt}
\right)
\!\!
\begin{pmatrix}
d\rho 	\\[0.5pt]
\sigma_1 	\\[0.5pt]
\sigma_2 	\\[0.5pt]
\sigma_3
\end{pmatrix}
\end{equation}
on the other (with, incidentally, the determinant of the transition matrix equal to $-4 \sqrt{\Delta}$), solving for $\chi_0$ in terms of $\chi_{+2}, \chi_{+1}, \chi_{-1}, \chi_{-2}$ is essentially the same as solving for $de_2$ in terms of $\sigma_1, \sigma_2, \sigma_3$ and $d\rho$. In practice, though, it is preferable to work with the $\chi$-frame. The remaining $m = \pm1$ components of equation \eqref{diff-GLT-eqs} together with $dx_{+2} = dz$ and the complex conjugate relation give us further the holomorphic differential frame $dz, du, d\bar{z}, d\bar{u}$ in terms of the frame $\chi_{+2}, \chi_{+1}, \chi_{-1}, \chi_{-2}$. Finally, we plug this result into the $j=2$ versions of the generalized Legendre transform metric and symplectic forms formulas \eqref{dimH=1-GLT-metric}--\eqref{dimH=1-GLT-sympl-forms}, with components given by the formulas \eqref{dimH=1-GLT-HA}. 

The computations involved at each step are quite long and laborious, but the end result is by comparison remarkably simple and elegant. To explain it, we need to first introduce a number of definitions. 

The first important observation is that  it is profitable to regard formally the $h_m$'s as the components of a spinor with spin 2. Notice that by construction they satisfy an alternating reality property, $\bar{h}_m = (-)^mh_{-m}$. From the middle equation \eqref{diff-GLT-eqs} it is clear that the natural way to introduce spherical normalization factors for them is by means of the \textit{dual} convention $h^2_m = (c^2_{m})^{-1} h^{\phantom{2}}_m$. If we define corresponding quadratic and cubic ``spherical invariants" by the formulas \eqref{rot-invs} with $x^2_m$ replaced by $h^2_m$ and take as above $N_2 = \sqrt{35/12}$, then in terms of the unnormalized components the formulas read
{\allowdisplaybreaks
\begin{equation} \label{r2(h)-r3(h)}
\begin{gathered}
r_2(h) = \frac{1}{2} \pt \big[ \pt |h_{+2}|^2 + 4 |h_{+1}|^2 + 6(h_0)^2 + 4 |h_{-1}|^2 + |h_{-2}|^2 \big] \\[0pt]
r_3(h) = 
\left|
\begin{array}{lll}
\! h_{+2} & h_{+1} & h_{\hp 0} \hspace{-3.5pt} \\
\! h_{+1} & h_{\hp 0} & h_{-1} \hspace{-3.5pt} \\
\! h_{\hp 0} & h_{-1} & h_{-2} \hspace{-3.5pt}
\end{array}
\right| \rlap{.}
\end{gathered}
\end{equation}
}%
These can be regarded as invariants only in a restricted sense. Namely, if we \textit{define}
\begin{equation} \label{h-ringed}
\mathring{h}^2_m = \sum_{m'=-2}^2 D^2_{mm'}(q^{-1}) \hp h^2_{m'}
\end{equation}
then we have indeed $r_2(h) = r_2(\mathring{h})$ and $r_3(h) = r_3(\mathring{h})$. 

Secondly, let us introduce the composite object
\begin{equation}
H^2_m = \frac{1}{2} \sqrt{\frac{7}{12}} 
\hspace{-1pt} 
\sum^2_{\shortstack{$\scriptstyle m_1,m_2=-2$ \\ $\scriptstyle m_1+m_2=m$}} 
\hspace{-3pt} 
\left(
\begin{array}{c|cc}
\! 2 & 2 & \!\!\! 2 \! \\
\! m & m_1 & \!\!\! m_2 \!
\end{array}
\right)
h^2_{m_1} h^2_{m_2}
\end{equation}
whose unnormalized components, assuming we spherically-normalize similarly to above by $H^2_m = (c^2_{m})^{-1} H^{\phantom{2}}_m$, are consequently given by the quadratic expressions
\begin{equation} 
\begin{aligned}
H_{+2} & = h_{+2}h_0 - h_{+1}h_{+1} \\[2pt]
H_{+1} & = \frac{1}{2} h_{+2}h_{-1} - \frac{1}{2} h_{+1}h_0 \\[-2pt]
H_{0\phantom{+}} & = \frac{1}{6}h_{+2}h_{-2} + \frac{1}{3} h_{+1}h_{-1} - \frac{1}{2} h_0h_0 \rlap{,}
\end{aligned}
\end{equation}
with the remaining ones easily determined by alternating conjugation. One can construct for it in the same way quadratic and cubic invariants, but these turn out to be, not surprisingly, composite as well:
\begin{equation} \label{r2(H)-r3(H)}
\begin{aligned}
r_2(H) & = \frac{1}{12} r_2(h)^2 \\
r_3(H) & = \frac{1}{4} r_3(h)^2 - \frac{1}{216} r_2(h)^3 \rlap{.}
\end{aligned}
\end{equation}
Due to the equivariance properties of the Clebsch-Gordan coefficients, it redefines covariantly at a redefinition \eqref{h-ringed}, that is
\begin{equation}
\mathring{H}^2_{m} = \sum_{m'=-2}^2 D^2_{mm'}(q^{-1}) H^2_{m'} \rlap{.}
\end{equation}

Thirdly, we introduce yet another frame by means of the linear transformation
\begin{equation} \label{chi-theta-frame}
\begin{pmatrix}
\chi_{+2} \\[1.5pt]
\chi_{+1} \\[1.5pt]
\chi_{-1}  \\[1.5pt]
\chi_{-2}
\end{pmatrix}
=
\left(
\hspace{1.1ex}
\begin{array}{@{}l@{\hskip 7pt}l@{\hskip 7pt}l@{\hskip 7pt}l@{}}
\mathring{H}'_{0} 			& \ 0 				& \ 0 				& \ 0 \\
\llap{$2$}\mathring{H}_{-1} 	& \mathring{H}'_0   	& \mathring{H}_{+2} 	& \ 0 \\
\ 0 						& \mathring{H}_{-2} 	& \mathring{H}'_{0} 	& \llap{$2$}\mathring{H}_{+1} \! \\
\ 0 						& \ 0 				& \ 0 			  	& \mathring{H}'_{0}
\end{array}
\right)
\!
\begin{pmatrix}
\theta_{+2} \\[1.5pt]
\theta_{+1} \\[1.5pt]
\theta_{-1}  \\[1.5pt]
\theta_{-2}
\end{pmatrix} 
\end{equation}
where, by definition,
\begin{equation}
\mathring{H}'_0 = \mathring{H}_0 - \frac{1}{6} r_2(\mathring{h}) \rlap{.}
\end{equation}
Note, incidentally, that $\mathring{H}'_0 \leq 0$. Indeed, by the definition of $r_2(\mathring{H})$ mirroring the first formula \eqref{r2(h)-r3(h)} one has $3 (\mathring{H}_0)^2 \leq r_2(\mathring{H})$, which then by the first equation \eqref{r2(H)-r3(H)}\,---\,and invariance\,---\,implies immediately that $(6\mathring{H}_0)^2 \leq r_2(\mathring{h})^2$. Since $r_2(\mathring{h}) > 0$, the claim follows. Observe also that the determinant of the transition matrix in \eqref{chi-theta-frame} is equal to $-r_3(\mathring{h}) \hp \mathring{h}_0(\mathring{H}'_0)^2$. 

Then the result of the computation outlined above can be stated as follows: 

\begin{proposition} \label{O4-Ansatz}
Consider the two successive linear transformations \eqref{quat-chi-frame} and \eqref{chi-theta-frame} of the quaternionic frame given by $\sigma_1,\sigma_2,\sigma_3$ and $d\rho$. In the resulting frame, the four-dimensional $j=2$ generalized Legendre transform hyperk\"ahler metric and symplectic forms assume the canonical form 
\begin{align}
g = \lambda \pt (\pt |\theta_{+1}|^2 + |\theta_{+2}|^2 \pt) \qquad\text{and}\qquad & \omega^1_m = \sum_{m'=-1}^1 D^1_{mm'}(q) \pt \mathring{\omega}^1_{m'} \\
& \text{with} \nonumber \\[2pt]
& \mathring{\omega}_{+1} =  \lambda \hp i \hspace{1.1ex} \theta_{+2} \wedge \bar{\theta}_{+1} \nonumber \\[0pt]
& \mathring{\omega}_{\hp 0} \hspace{6.1pt} = \lambda \hp i \pt (\hp \theta_{+1} \wedge \bar{\theta}_{+1} - \theta_{+2} \wedge \bar{\theta}_{+2}) \nonumber \\[0pt]
& \mathring{\omega}_{-1} = \lambda \hp i \hspace{1.1ex} \bar{\theta}_{+2} \wedge \theta_{+1} \nonumber
\end{align}
where $\lambda = r_3(\mathring{h}) \mathring{H}'_{0}$ and we normalize $\omega^1_m = c^1_m \omega^{\phantom{1}}_m$. 
\end{proposition}


The parallel with the form \eqref{HKstr-j=1} in the $j=1$ case is evident and demonstrates that quaternionic parametrization of the Majorana polynomial coefficients opens indeed a path to generalizing the Gibbons-Hawking Ansatz. 

The hyperk\"ahler metric and symplectic forms are completely determined modulo simple algebraic operations by the five functions $(\mathring{h}_m)_{m=-2,\dots,2}$, closely related to the functions $(h_m)_{m=-2,\dots,2}$, which in turn can be extracted from a single holomorphic function through a contour-integral of the type \eqref{h-bar-oint}. Note that, by comparison, the generalized Legendre transform formulas \eqref{dimH=1-GLT-HA} with $j=2$ depend in addition to these also on $h_{-3}$. The contour-integral formula gives the $h_m$ functions naturally in terms of the variables $x_m$, from which passing to the quaternionic variables in which the metric and symplectic forms are expressed is usually straightforward. In the generalized Legendre transform formulas \eqref{dimH=1-GLT-metric}--\eqref{dimH=1-GLT-sympl-forms}, by contrast, the metric and symplectic forms are expressed in a holomorphic coordinate frame, and the dependence of the variables $x_m$ on the holomorphic coordinates is virtually always implicit. In either case one still has to impose of course the generalized Legendre transform constraint $L_{x_0} = 0$. 

In practice, therefore, one needs to compute two things: the first derivative $L_{x_0}$, which determines the generalized Legendre transform constraint, and the second derivatives $L_{x_0x_m}$ with $m=-2,\dots,2$, which then determine the hyperk\"ahler metric and symplectic forms.

\section{The $\mathcal{O}(4)$ spectral curve} \label{sec:O4-curve}

In this section we continue our generic exploration of the generalized Legendre transform construction with $j=2$ with a detailed study of the structure of a certain elliptic curve associated to the characteristic  $\mathcal{O}(4)$ section which plays an important role in applications. 

\subsection{} 

In the total space of the bundle $\pi^*\mathcal{O}(4) \rightarrow Z$ with global section $x$ consider the quartic plane curve given by
\begin{equation} \label{O4-sc-Maj}
\eta^2 = x_N(\zeta)
\end{equation}
where $\eta$ is the complex coordinate on the fiber and
\begin{equation}
x_N(\zeta) = \sum_{m=-2}^{2} x_m \hp \zeta^{2-m}
\end{equation}
with coefficients satisfying the alternating reality condition $\bar{x}_m = (-)^m x_{-m}$ is the usual $j=2$ Majorana polynomial representing the section $x$ in the local trivialization over the open set $N \subset Z$. We call such a curve an $\mathcal{O}(4)$ \textit{spectral curve}, and label its defining representation as \textit{Majorana normal form}. 

By construction, this curve has two involutions, one holomorphic and one anti-ho\-lo\-mor\-phic:
\begin{equation} \label{invols}
\begin{aligned}
& \text{elliptic involution:} && (\zeta,\eta) \mapsto (\zeta, - \eta) \\
& \text{real structure:} && (\zeta,\eta) \mapsto (-1/\bar{\zeta},-\bar{\eta}/\bar{\zeta}^2) \rlap{.}
\end{aligned}
\end{equation}
One should think of $\eta$ as the complex square root of the quartic Majorana polynomial, and, indeed, the introduction of the elliptic curve can be justified by the need to cast the use of this square root on a rigorous basis. The elliptic involution encapsulates the branch ambiguity data for the square root, and the second rule essentially implies that under antipodal conjugation (which, recall, we define here as the substitution $\zeta \mapsto \zeta^c$ followed by complex conjugation) the square root changes branches, that is, $\sqrt{x_N} \mapsto - \sqrt{x_S}$. 

One also has a natural nowhere-vanishing holomorphic 1-form on the curve, which is to say, an abelian differential of the first kind, given by
\begin{equation}
\theta = \frac{d\zeta}{2 \hp \eta} \rlap{.}
\end{equation}

Let $a_1,\dots,a_4$ denote the roots of the Majorana quartic polynomial. The reality properties of the polynomial coefficients imply that they come in antipo\-dally-conjugated pairs or, in other words, the two sets $\{ a_1, \dots, a_4 \}$ and $\{a^c_1, \dots, a^c_4\}$ are the same. In terms of them we have 
\begin{equation}
\eta^2 = x_{-2} \prod_{i=1}^4 (\zeta - a_i) \rlap{.}
\end{equation}
We can cast the curve in the \textit{Weierstrass normal form}
\begin{equation} \label{W.n.f.}
Y^2 = X^3 - g_2 X - g_3
\end{equation}
by taking for example the root $a_4$ to $\infty$ by means of the birational transformation
\begin{align}
X = x_{-2} \bigg( S_4 - \frac{P_4}{\zeta-a_4} \bigg) \quad \text{with} & \quad \cramped{S_4 = \frac{1}{3} \sum_{1\leq i < j \leq 3} (a_i-a_4)(a_j-a_4)} \\
& \quad \cramped{P_4 = \prod_{1 \leq i \leq 3} (a_i-a_4)} \rlap{.} \nonumber
\end{align}
and then choosing an appropriate transformation for $\eta$. The remaining Majorana roots are mapped to the Weierstrass roots $e_1, e_2, e_3$. Note in particular that with our ordering convention for the Weierstrass roots, $a^c_4$, the root antipodally conjugated to $a_4$, is mapped to $e_2$. The Weierstrass coefficients $g_2$ and $g_3$ which arise here are precisely those defined in equation \eqref{Weierstrass-coeffs}\,---\,and this brings us again, through a different route, to the same cubic Weierstrass polynomial that we have encountered when we introduced a quaternionic parametrization for the quartic Majorana coefficients $x_m$. 

In the new coordinates the holomorphic 1-form takes the form
\begin{equation} \label{M-W}
\theta = \frac{d\zeta}{2 \hp \eta} = \frac{dX}{2 Y} \rlap{.}
\end{equation}

Finally, by means of the Weierstrass $\wp$-function and its derivative we pass to a complex torus representation. Taking $X = \wp(u;4g_2,4g_3)$ and $2Y = \wp'(u;4g_2,4g_3)$ gives further
\begin{equation} \label{1-f-MWJ}
\theta = \frac{d\zeta}{2 \hp \eta} = \frac{dX}{2 Y} = du \rlap{.}
\end{equation}
The reality of the Weierstrass roots implies that the toric lattice $\Lambda$ is orthogonal. If $\omega_1$, $\omega_2$, $\omega_3$ with $\omega_1+\omega_2+\omega_3=0$ are the Weierstrass \textit{half}-periods, then our ordering of the Weierstrass roots entails that $\omega_1$ is real and $\omega_3$ purely imaginary. Following the usual conventions in the literature, in what follows we will often denote these as $\omega$ and $\omega'$, respectively. 

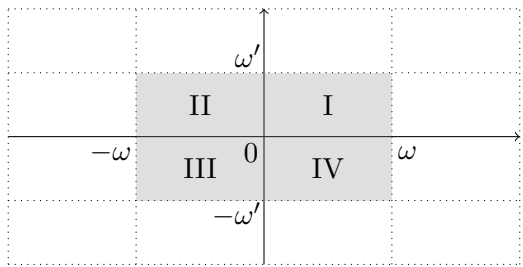
\begin{figure}[h]
\begin{tikzpicture}[scale=0.85]

\draw[dotted] (-4,2) -- (4,2);
\draw[dotted] (-4,1) -- (4,1);
\draw[->] (-4,0) -- (4,0);
\draw[dotted] (-4,-1) -- (4,-1);
\draw[dotted] (-4,-2) -- (4,-2);

\draw[dotted] (-4,-2) -- (-4,2);
\draw[dotted] (-2,-2) -- (-2,2);
\draw[->] (0,-2) -- (0,2);
\draw[dotted] (2,-2) -- (2,2);
\draw[dotted] (4,-2) -- (4,2);

\path[fill=gray, nearly transparent] (-2,1) -- (2,1) -- (2,-1) -- (-2,-1) -- (-2,1);

\draw (-2.22,-0.25) node {$\llap{$-$}\omega$};
\draw (-0.20,-0.25) node {$0$};
\draw ( 2.25,-0.25) node {$\omega$};

\draw (-0.25,1.27) node {$\omega'$};
\draw (-0.25,-1.22) node {$\llap{$-$}\omega'$};

\draw (1,0.5) node {I};
\draw (-1,0.5) node {II};
\draw (-1,-0.5) node {III};
\draw (1,-0.5) node {IV};

\end{tikzpicture}
\caption{A choice of fundamental domain for the toric lattice $\Lambda$, divided in quadrants. We assume that the domain has closed upper and right boundaries and open lower and left ones.} \label{lattice-fund-dom}
\end{figure}

\subsection{} \label{ssec:ac-torus}

Let us try to understand now how the two involutions \eqref{invols} act on the complex torus $\mathbb{C}/\Lambda$. Given an arbitrary point on the curve with Majorana normal coordinates $(\zeta,\eta)$, let $(X_{\zeta},Y_{\zeta})$ denote its image in Weierstrass normal coordinates, and $u_{\zeta}$ its image on the complex torus. 

The elliptic involution induces, modulo a lattice shift, a reflection with respect to the origin of the complex plane: 
\begin{equation} \label{ell-inv}
u_{\zeta} \mapsto  - u_{\zeta} \mod \Lambda \rlap{.}
\end{equation}
If we require $u_{\zeta}$ to be in the fundamental domain of the lattice as defined in Figure \ref{lattice-fund-dom}, then so is $-u_{\zeta}$, and consequently we can dispense in \eqref{ell-inv} with lattice shifts altogether. Note that in this case the map $u_{\zeta} \mapsto - u_{\zeta}$ interchanges diagonally-opposite quadrants. 

To work out the action of antipodal conjugation, we begin by noticing the following useful identity
\begin{equation} \label{ambig_sign_fix}
\frac{Y_{\zeta}}{X_{\zeta} - X_{\zeta'}} = \eta  \bigg( \frac{1}{\zeta-\zeta'} - \frac{1}{\zeta-a_4} \bigg)
\end{equation}
valid for any $\zeta' \in \mathbb{C}$ different from $\zeta$. This can be checked by squaring each side up and then imposing both the Majorana and Weierstrass forms of the elliptic equation. A possible sign ambiguity is fixed by the requirement of consistency in the limit $\zeta' \rightarrow \zeta$ with the second equation \eqref{M-W}. 

Taking in particular $\zeta' = a^c_4$, we obtain
\begin{equation}
\frac{Y_{\zeta}}{X_{\zeta} - e_2} = - \eta \pt \frac{1+|a_4|^2}{(\zeta - a_4)(1 + \bar{a}_4\zeta)} \rlap{.}
\end{equation}
The right-hand side of this equation is manifestly invariant under antipodal conjugation and so consequently
\begin{equation}
\frac{Y_{\zeta}}{X_{\zeta} - e_2} = \frac{\bar{Y}_{\zeta^c}}{\bar{X}_{\zeta^c} - e_2} \rlap{.}
\end{equation}
This can be equivalently rephrased as the \textit{collinearity} condition
\begin{equation}
\left|
\begin{array}{cll}
1 & X_{\zeta} & Y_{\zeta} \\[2pt]
1 & \bar{X}_{\zeta^c}\!\! & \bar{Y}_{\zeta^c} \!\! \\
1 & \pt e_2 & \pt 0
\end{array}
\right| = 0
\end{equation}
which then by way of the addition law of the Weierstrass elliptic functions allows us to conclude that
\begin{equation}
u_{\zeta} + \bar{u}_{\zeta^c} + \omega_2 = 0 \mod \Lambda \rlap{.}
\end{equation}
The mod $\Lambda$ ambiguity can be fixed by imposing some restricting condition, typically a choice of fundamental domain for the lattice, in which case we can write
\begin{equation} \label{bar-u-ac}
u_{\zeta} + \bar{u}_{\zeta^c} = n^{\phantom{'}}_{\zeta} \omega + n'_{\zeta}\omega'
\end{equation}
for some, very importantly, \textit{odd} integers $n^{\phantom{'}}_{\zeta}$ and $n'_{\zeta}$. For instance, if we require that both $u_{\zeta}$ and $u_{\zeta^c}$ be in the fundamental rectangle of the toric lattice as defined in Figure \ref{lattice-fund-dom}, then we find that the shift term is entirely determined by the particular \textit{quadrant} of the fundamental domain to which, say, $u_{\zeta}$ belongs, as follows:
\begin{equation*}
\begin{tabular}{c|c|c}
$u_{\zeta} \in $ quadrant & \hspace{1ex} $n_{\zeta}$ \hspace{1ex} & \hspace{1ex} $n'_{\zeta}$ \hspace{1ex} \\[2pt] \hline
I & $1$ & $1$ \\ \hline
II & $\llap{$-$}1$ & $1$ \\ \hline
III & $\llap{$-$}1$ & $\llap{$-$}1$ \\ \hline
IV & $1$ & $\llap{$-$}1$ 
\end{tabular}
\end{equation*}
Note that, quite generally, $n^{\phantom{'}}_{\zeta^c} = n^{\phantom{'}}_{\zeta}$ and $n'_{\zeta^c} = - n'_{\zeta}$, and so the map $u_{\zeta} \mapsto u_{\zeta^c}$ interchanges in this case vertically-adjacent quadrants. A shift of the fundamental domain by $\lambda \pt \omega$ in the real direction and $\lambda'\omega'$ in the imaginary one for some $\lambda,\lambda' \in \mathbb{Z}$ entails a shift $n^{\phantom{'}}_{\zeta} \mapsto n^{\phantom{'}}_{\zeta} + 2\hp \lambda$ while leaving $n'_{\zeta}$ unchanged. For reasons of generality, in the remainder of the paper we will assume a generic (yet unspecified) such choice of fundamental domain. 

Observe that we can define an ``enhanced" version of $u_{\zeta}$ by $\mathfrak{u}_{\zeta} = u_{\zeta} - (n^{\phantom{'}}_{\zeta}/2)\hp \omega - (n'_{\zeta}/2)\hp \omega'$. This modified variable satisfies the simpler antipodal conjugation property $\bar{\mathfrak{u}}_{\zeta^c} = - \mathfrak{u}_{\zeta}$ and, moreover, its real part is independent of the particular choice of fundamental domain and takes values the interval $(-\omega/2,\omega/2]$. Let furthermore 
\begin{equation}
\begin{aligned}
\upsilon_{\zeta} & = \mathfrak{u}_{\zeta} + \bar{\mathfrak{u}}_\zeta = u_{\zeta} - u_{\zeta^c} - n'_{\zeta}\omega' \\
\upsilon'_{\zeta} & = \mathfrak{u}_{\zeta} - \bar{\mathfrak{u}}_\zeta = u_{\zeta} + u_{\zeta^c} - n_{\zeta}\omega \rlap{.}
\end{aligned}
\end{equation}
Thus defined, $\upsilon_{\zeta}$ is real and $\upsilon'_{\zeta}$ purely imaginary. As the $g_2$, $g_3$ parameters are real, all standard Weierstrass functions have definite reality properties at these points which one can easily infer based on their parity. Hence, for example, since the Weierstrass $\wp$-function and its derivative are even respectively odd functions, it follows that both $\wp(\upsilon_{\zeta})$ and $\wp(\upsilon'_{\zeta})$ are real, while $\wp'(\upsilon_{\zeta})$ is real and $\wp'(\upsilon'_{\zeta})$ purely imaginary. Based on the elliptic formulas \eqref{assoc-sigma-propr} and similar arguments relying on the fact that the Weierstrass sigma-function is odd and all three associated sigma-functions are even, we can argue moreover that for all values of $\zeta$ we have
\begin{equation}
\wp(\upsilon'_{\zeta}) \leq e_3 < e_2 < e_1 \leq \wp(\upsilon_{\zeta}) \rlap{.}
\end{equation}

\subsection{}

In the remainder of this section we discuss some further properties of $\mathcal{O}(4)$ spectral curves which, while interesting in themselves, have no immediate bearing on the considerations which follow. The pragmatic reader is urged to skip ahead to the next section. 

Taking now in the relation \eqref{ambig_sign_fix} $\zeta=0$ and $\zeta'=\infty$ and then vice versa, and noticing that we have
\begin{equation}
X_{0} = \frac{x_0}{3} + \frac{x_{+1}}{a_4} + \frac{2 \hp x_{+2}}{a_4^2}
\quad\text{and}\quad
X_{\infty} = \frac{x_0}{3} + x_{-1}a_4 + 2 \hp x_{-2}a_4^2 \rlap{,}
\end{equation}
we can derive two more collinearity properties, namely,
\begin{equation}
\left|
\begin{array}{cll}
1 & X_{0} & \, Y_{0} \\[2pt]
1 & X_{\infty} & \, Y_{\infty} \! \\
1 & \pt x_{+} & \, \llap{$-$} \pt y_{+}
\end{array}
\right| = 0
\qquad\text{and}\qquad
\left|
\begin{array}{cll}
1 & X_{0} & \, Y_{0} \\[2pt]
1 & X_{\infty} & \, \llap{$-$} Y_{\infty} \! \\
1 & \pt x_{-} & \, \pt y_{-}
\end{array}
\right| = 0
\end{equation}
where, by definition,
\begin{equation}
\begin{aligned}
x_{\pm} & = \frac{x_0}{3} \pm 2\sqrt{x_{-2}}\sqrt{x_{+2}} \\
y_{\pm} & = x_{+1} \sqrt{x_{-2}} \pm x_{-1} \sqrt{x_{+2}} \rlap{.}
\end{aligned}
\end{equation}
The complex square roots $\sqrt{x_{+2}}$ and $\sqrt{x_{-2}}$ arise from $\eta$ at $\zeta =0$ and $\infty$, and their relative branch choice ambiguity is fixed by requiring them to be complex conjugates. So then what the collinearity relations show is that 
\begin{equation}
\begin{aligned}
x_+ & = \hspace{2.6pt} \wp(u_{\infty}+u_0) 	\qquad &  x_- & = \hspace{2.6pt} \wp(u_{\infty}-u_0) \\
2y_+ & = \wp'(u_{\infty}+u_0) 	& 2y_- & = \wp'(u_{\infty}-u_0) \rlap{.}
\end{aligned}
\end{equation}

From their definitions it is clear that both $x_-$ and $x_+$ are real, while $y_-$ is real and $y_+$ purely imaginary. Moreover, we have the following inequalities:
\begin{equation}
\begin{gathered}
e_3 < x_- < e_2 < x_+ < e_1 \\
e_3 < - x_+ - x_- < e_1 \rlap{.}
\end{gathered}
\end{equation}
To prove them, we rely on the representation \eqref{wv-Majorana} of Majorana polynomials. By expressing the Majorana coefficients and Weierstrass roots in terms of the complex variables $w_i,v_i$, we can derive the following set of formulas:
{\allowdisplaybreaks
\begin{equation}
\begin{aligned}
x_{\pm} - e_1 & = - (|v_1w_2| \mp |v_2w_1|)^2 \\
x_{\pm} - e_2 & = - (\sqrt{v_1w_1\bar{v}_2\bar{w_2}} \mp \sqrt{\bar{v}_1\bar{w}_1v_2w_2})^2 \\
x_{\pm} - e_3 & = \phantom{+} (|v_1v_2| \pm |w_1w_2|)^2
\end{aligned}
\end{equation}
}%
and
{\allowdisplaybreaks
\begin{equation}
\begin{aligned}
x_+ + x_- + e_1 & = \phantom{+} |v_1\bar{v}_2 - w_2\bar{w}_1|^2 \\
x_+ + x_- + e_2 & = \phantom{+} (|v_1|^2 - |w_1|^2)(|v_2|^2 - |w_2|^2) \\
x_+ + x_- + e_3 & = - |v_1\bar{w}_2 + v_2\bar{w}_1|^2 \rlap{,}
\end{aligned}
\end{equation}
}%
from which then the inequalities follow promptly. 

The pairs $(x_-,y_-)$ and $(x_+,y_+)$ are points on the Weierstrass curve, as one may easily convince oneself directly by verifying that they both satisfy indeed the Weierstrass equation \eqref{W.n.f.}. Turning this around, we can solve the two resulting equations for the Weierstrass coefficients, to obtain the alternative expressions
\begin{equation}
\begin{aligned}
g_2 & = x_+^2 + x_+x_- + x_-^2 - \frac{y_+^2 - y_-^2}{x_+ - x_-} \\
g_3 & = -x_+x_-(x_+ + x_-) + \frac{x_-y_+^2 - x_+y_-^2}{x_+ - x_-} \rlap{.}
\end{aligned}
\end{equation}
These expressions allow us then to immediately see that the Weierstrass cubic polynomial associated to the Majorana quartic polynomial can be cast in the following  determinantal form:
\begin{equation}
X^3 - g_2X - g_3 =
\begin{vmatrix}
X-x_+ & \displaystyle \frac{-iy_+}{\sqrt{x_+ - x_-}} & 0 \\
\displaystyle \frac{-iy_+}{\sqrt{x_+ - x_-}} & X + x_+ + x_- & \displaystyle \frac{y_-}{\sqrt{x_+ - x_-}} \\
0 & \displaystyle \frac{y_-}{\sqrt{x_+ - x_-}} & X - x_-
\end{vmatrix} 
\end{equation}
with each entry, apart from the variable $X$, real. Note that the formula remains valid if we take in it, either separately or together, $y_+ \mapsto - y_+$ and $y_- \mapsto - y_-$.

\section{ALE and ALF gravitational instantons of type $A_k$}

We shift now our focus to applications, and our main goal for the remainder of the paper will be the explicit construction of gravitational instantons of type $D_k$. However, before we tackle that, we find it instructive to review Hitchin's construction of ALE gravitational instantons of type $A_k$ \cite{MR520463}, which will serve as a model for our treatment of the $D_k$ case. Our review, though, will have a slightly skewed objective compared to the original: we will use the approach of \cite{MR520463} in order to derive the input data for the Legendre transform construction\,---\,which consists, mostly, of a holomorphic function\,---, and then use this rather than Hitchin's original method to construct the hyperk\"ahler metric and symplectic forms. The result will evidently be the same, but the variation in strategy is really a practice run for the $D_k$ case, where we will follow the same blueprint. 

\subsection{}

The starting point of Hitchin's approach is the minimal resolution of a complex algebraic surface with a Kleinian singularity of type $A_k$. A standard construction going back to F. Klein shows that the singular quotients $\mathbb{C}^2/\Gamma$ can be embedded in $\mathbb{C}^3$ for any finite subgroup $\Gamma$ of $SU(2)$. In particular, if $\Gamma$ is the cyclic group of order $k+1$, then the resulting singularity is known as a Kleinian (or du Val, or rational double point) singularity of type $A_k$, and may be described by the following algebraic equation in $\mathbb{C}^3$: 
\begin{equation}
y^2 + z^2 = x^{k+1} \rlap{.}
\end{equation}
In terms of a set of $k+1$ constant complex parameters $t_l$ satisfying the constraint $\sum_{l=0}^kt_l=0$, its universal deformation has the form
\begin{equation} \label{Ak-def}
y^2 + z^2 = \prod_{l=0}^k (x + t_l) \rlap{.}
\end{equation}
This non-singular complex surface comes equipped with a natural nowhere-vanishing symplectic 2-form given by
\begin{equation}
\varpi = \frac{2 \pt dy \wedge dz}{f_x} = \frac{2 \pt dz \wedge dx}{f_y} = \frac{2 \pt dx \wedge dy}{f_z}
\end{equation}
where $f(x,y,z)$ is the function whose vanishing locus defines the surface, and the factor 2 has been chosen for ulterior convenience. The first formula is valid on the patch defined by the condition $f_x \neq 0$, where the variables $y$ and $z$ give a good system of coordinates on the surface, the second one on the patch defined by $f_y \neq 0$, and so on. 

Using a hyperk\"ahler quotient construction, Kronheimer has shown in \cite{MR992334} that one can define on the desingularized surface an ALE hyperk\"ahler metric. So then we may ask ourselves the following question\,---\,and this is really the key to the the twistor approach: the above complex holomorphic description of the algebraic surface is holomorphic with respect to \textit{which} particular hyperk\"ahler complex structure, out of the many we have at our disposal? The answer that Hitchin essentially proposed was that the surface can be described by this same algebraic equation in almost all hyperk\"ahler complex structures. In other words, one can think of the variables $x,y,z$ and also of the parameters $t_l$ in some appropriate sense to be defined as \textit{sections} of certain holomorphic line bundles over the twistor space of the complex surface viewed as a hyperk\"ahler space. We take then in particular $x = x_{NS}(\zeta)$ and, redenoting the parameters $t_l \equiv - x_l$, $x_l=x_{l,NS}(\zeta)$ to be the tropical parts of real (with respect to antipodal conjugation) holomorphic global sections of the pull-back bundle $\pi^*\mathcal{O}(2)$.

To determine the structure of the remaining variables $y$ and $z$ we begin by observing that if we introduce in their place a new set of complex variables $\xi_{N,S}  = y \pm i z$ then the equation \eqref{Ak-def} factorizes as follows
\begin{equation} \label{xiNS-x}
\xi_N \xi_S = \prod_{l=0}^k (x - x_l)
\end{equation}
and the holomorphic 2-form becomes
\begin{equation} \label{om-xi-x}
\varpi = -i \pt d \pt ( \hp \ln \xi_N ) \wedge \hspace{-0.5pt} dx = i \pt d \pt ( \hp \ln \xi_S ) \wedge \hspace{-0.5pt} dx 
\end{equation}
with each expression valid on some appropriately defined patch. Switching further to the variables $p_N$ and $p_S$ defined by
\begin{equation} \label{xisps}
\begin{aligned}
\xi_N & = e^{+ i \hp p_N} \\
\xi_S \pt & = e^{ - i \hp p_S}
\end{aligned}
\end{equation}
recasts the equation in the form
\begin{equation} \label{exp(pN-pS)}
e^{i\hp (p_N-p_S)} = \prod_{l=0}^k (x - x_l) \rlap{.}
\end{equation}
Notice that we can write this equivalently as
\begin{equation} \label{pN-pS-j=1}
p_N - p_S = i \frac{\partial H(x)}{\partial x}
\end{equation}
where, by definition,
\begin{equation}
\frac{\partial H(x)}{\partial x} = -  \ln \prod_{l=0}^k (x - x_l) \rlap{.}
\end{equation}
Integrating gives us, up to an integration constant,
\begin{equation} \label{ALE-H-Ak}
H(x) = - \sum_{l=0}^k [\pt  (x - x_l)\ln(x - x_l) - (x - x_l)] \rlap{.}
\end{equation}
In the current variables the holomorphic 2-form takes the form
\begin{equation} \label{twist-Darboux}
\varpi = \frac{dp_N \wedge dx_N}{\zeta} = \zeta \pt dp_S \wedge dx_S
\end{equation}
where, recalling the definition \eqref{arct-trop-anta}, we have replaced successively the tropical form of $x$ with its arctic respectively antarctic ones. 

Put in this way, this then conspicuously suggests that we identify the nowhere-vanishing holomorphic 2-form with the fiberwise-supported holomorphic twistor space 2-form as follows
\begin{equation}
\varpi = \omega_{NS}(\zeta) 
\end{equation}
and at the same time the $p$ and $x$-variables with the similarly named Legendre transform ones. In this case, the equation \eqref{pN-pS-j=1} is precisely the equation \eqref{pN-pS} in real dimension four with $j=1$, and $H(x)$ is the holomorphic function which in the the Legendre transform construction characterizes the hyperk\"ahler structure. 

\subsection{}

In the Legendre transform approach critical regularity conditions require $p_N$ and $p_S$ to be well-defined at $\zeta=0$ and $\infty$, respectively\,---\,and if they are to be satisfied we have to have $\xi_N$ and $\xi_S$ obey the same conditions. Moreover, $\xi_N$ and $\xi_S$ are expected to be interchanged by antipodal conjugation. Let us show now that such solutions of the equation \eqref{xiNS-x} exist. These conditions are in fact so stringent that they determine the $\xi$-variables almost completely.

As the tropical component of a real section of $\pi^*\mathcal{O}(2)$, $x_l-x$ admits, by virtue of the same argument used to justify 
the representation \eqref{Majorana-roots} of Majorana polynomials, a factorization of the type
\begin{equation}
x_l - x = \rho_l \frac{(\zeta-a_l)(1+\bar{a}_l\zeta)}{(1+|a_l|^2)\hp \zeta}
\end{equation}
with $\rho_l > 0$. So if we take
\begin{equation} \label{xi-expl}
\begin{aligned}
\xi_N & = \prod_{l=0}^k \sqrt{\rho_l} \pt \frac{a_l-\zeta}{\sqrt{1+|a_l|^2}} \\
\xi_S \pt\hp & = \prod_{l=0}^k \sqrt{\rho_l} \pt \frac{\bar{a}_l+\tilde{\zeta}}{\sqrt{1+|a_l|^2}} \\
\end{aligned}
\end{equation}
where $\tilde{\zeta} = 1/\zeta$, then these are 1.~manifestly interchanged by antipodal conjugation, 2.~well-defined at $\zeta=0$ respectively $\zeta=\infty$ and 3.~they automatically satisfy the equation \eqref{xiNS-x}. This proves existence. As far as uniqueness is concerned, observe that $\xi_N$ and $\xi_S$ are defined in principle only up to (opposite) constant complex phases. 

\subsection{} 

It is known, see \textit{e.g.} \cite{MR877637}, that the generalized Legendre transform holomorphic $H$-functions describing the cyclic ALF gravitational instantons can be obtained from their ALE counterparts \eqref{ALE-H-Ak} very simply through the addition of a quadratic term, namely
\begin{equation} \label{H^ALF-j=1}
H^{ALF}(x) = - \alpha \hp x^2 + H^{ALE}(x)
\end{equation}
where $\alpha$ is a positive real constant inversely related to the so-called ``mass parameter". In the ALF case one is then led to consider the \textit{deformed} equation
\begin{equation}
\xi_N\xi_S = e^{2 \hp \alpha \hp x} \prod_{l=0}^k (x-x_l) \rlap{.}
\end{equation}
Notice that the formulas \eqref{om-xi-x} still give without any modification a nowhere-vanishing holomorphic 2-form on the deformed surface. If we then proceed to define $p_N$ and $p_S$ by the same formulas, we obtain the same equation \eqref{pN-pS-j=1} but now with the modified $H$-function \eqref{H^ALF-j=1}, while the holomorphic 2-form assumes in the end the same twisted Darboux form \eqref{twist-Darboux}. To work out how the formulas \eqref{xi-expl} are modified, we write $x$ in Majorana form and then split it into the sum of two antipodally-conjugated terms as follows:
\begin{equation} \label{x-split-Ak}
x = \underbracket[0.4pt][4pt]{ \frac{x_+}{\zeta\ } + \frac{x_0}{2} }_{\text{sing.\ at $\zeta = 0$}} + \underbracket[0.4pt][6.1pt]{ \frac{x_0}{2} + x_-\zeta }_{\text{sing.\ at $\zeta = \infty$}} \rlap{\!\!.}
\end{equation}
This gives us promptly
\begin{equation}
\begin{aligned}
\xi^{ALF}_N & = e^{2 \alpha \hp \big(\frac{x_0}{2} + x_-\zeta \big)} \xi^{ALE}_N \\
\xi^{ALF}_S & = e^{2 \alpha \hp \big( \frac{x_+}{\zeta\ } \pt + \pt \frac{x_0}{2} \big)} \pt\xi^{ALE}_S 
\end{aligned}
\end{equation}
satisfying all the required properties.

\subsection{} As far as the effective construction of the metric is concerned, the essential outcome of this argument is the $H$-function, which constitutes the main input in the Legendre transform machinery. From this, one can compute the contour integral $L$ explicitly by choosing a contour in such a way as to obtain a real answer, and then its second derivatives. The end result can be cast in Gibbons-Hawking form with the well-known multi-center $\mathbb{R}^3$-harmonic potential 
\begin{equation}
V = \alpha + \sum_{l=0}^k \frac{1}{|\vec{r} - \vec{r}_l|} \rlap{.}
\end{equation}

\subsection{The period matrix} 

Let us return now to the ALE case to compute the period matrix. The exceptional divisor of the minimal resolution $M$ of the Kleinian singularity $\mathbb{C}^2/\Gamma$ with $\Gamma$ the cyclic group of order $k+1$ on which the metric is defined consists of $k$  rational curves whose homology classes form a basis in $H_2(M,\mathbb{Z})$, with intersection matrix given by minus the Cartan matrix of the simply-laced Lie algebra $A_k$. We can pair these exceptional 2-cycles with the cohomology classes in $H^2(M,\mathbb{R})$ corresponding to the three standard hyperk\"ahler symplectic forms to obtain the so-called \textit{period matrix}. In the $A_k$ case this has been computed by Hitchin in \cite{MR520463}. Here, however, we will use a different method, one which we can generalize easier to the $D_k$ case. 

The deformed algebraic equation \eqref{Ak-def} can be regarded as describing a fibration of complex quadrics over the complex $x$-plane, with the fibers degenerating at the points $x=-t_l=x_l$. For any pair of distinct indices $l_1,l_2 \in \{0,\dots,k\}$ let us consider inside this complex surface a real 2-cycle $\gamma_{\hp l_1l_2}$ in the form of a circle fibration over a 1-dimensional contour in the $x$-plane connecting the points $x_{l_1}$ and $x_{l_2}$, with the circle degenerating to a point at the ends. In terms of the $x,p_N,p_S$ variables we parametrize the 2-cycle as follows:
\begin{equation}
\begin{aligned}
& \text{$x  = x(\lambda)$ such that $x(0)=x_{l_2}$ and $x(1)=x_{l_1}$} \\
& p_N = \phi + p_N^0(\lambda) \\
& p_S\hspace{1.9pt} = \phi + \hp p_{\hp S}^0\hp(\lambda)
\end{aligned}
\end{equation}
with $\lambda \in [0,1]$ and $\phi \in [0,2\pi]$ denoting the base and $S^1$-fiber parameters, respectively. The structure of this parametrization does not change under a change of complex structure. This fact relies on two important checks: First, the dependence of $p_N$ and $p_S$ on $\phi$ is consistent with antipodal conjugation. Second, the lack of dependence of $x$ on $\phi$ is consistent with the cancellation of $\phi$ in the exponent on the left-hand side of equation \eqref{exp(pN-pS)}. 

Using the complex Darboux form of the holomorphic twistor space 2-form and marking the dependence on $\zeta$ explicitly, we have then successively 
{\allowdisplaybreaks
\begin{align}
\int_{\gamma_{\hp l_1l_2}} \! \omega_N(\zeta) 
& = \int_{\gamma_{\hp l_1l_2}} dp_N(\zeta) \wedge dx_N(\zeta) \\[2pt] 
& = \int_{\gamma_{\hp l_1l_2}} d\phi \wedge dx_N(\zeta) \nonumber \\[-2pt] 
& = \int_0^{2\pi} \!\! d\phi \int_{x_{l_2N}(\zeta)}^{\mathstrut  x_{l_1N}(\zeta)} dx_N(\zeta) \nonumber \\[5pt]
& = 2\pi \hp [\hp x_{l_1N}(\zeta) - x_{l_2N}(\zeta)] \rlap{.} \nonumber
\end{align}
}%
The second line follows from the fact that, when restricted to $\gamma_{\hp l_1l_2}$, both $x_N$ and $p_N^0$ depend only on $\lambda$ and thus give a vanishing $d\lambda \wedge d\lambda$ term. Since $\zeta$ is arbitrary and since, as we have argued above, the definition of $\gamma_{\hp l_1l_2}$ does not depend on the choice complex structure and thus of $\zeta$, we infer that
\begin{align}
\int_{\gamma_{l_1l_2}} \vec{\omega} = 2\hp \pi \hp (\vec{r}_{l_1} - \vec{r}_{l_2}) \rlap{.}
\end{align}
The homology classes of the 2-cycles $\gamma_{\hp l_1l_2}$ are not all linearly independent, but we can always choose from among them a basis $[\gamma_i\hp]$ in $H_2(M,\mathbb{Z})$ with elements relabeled by $i = 1,\dots,k$ such that
\begin{align}
\int_{\gamma_{i}} \vec{\omega} = 2\hp \pi \hp (\vec{r}_{i-1} - \vec{r}_i) \rlap{.}
\end{align}

\section{ALE and ALF gravitational instantons of type $D_k$}

We are now finally ready to address directly the problem of constructing gravitational instantons of type $D_k$. The strategy of our approach will essentially mirror the one we laid out in the $A_k$ case. First, we will show, using crucial ideas from Cherkis-Kapustin \cite{MR1693628,MR1700937} and Cherkis-Hitchin \cite{MR2177322} and the elliptic apparatus that we have developed in section \ref{sec:O4-curve}, how one can derive starting from the minimal resolution of a Kleinian $D_k$ singularity the generalized Legendre transform holomorphic $H$-function characterizing the instanton. The relevant generalized Legendre transform construction is the $\mathcal{O}(4)$ (or $j=2$) one, and the argument yields also naturally the characteristic constraint in the form of a divisor class condition on the $\mathcal{O}(4)$ spectral curve. To determine the hyperk\"ahler metric and symplectic forms we will rely on the $\mathcal{O}(4)$ Ansatz we have formulated in section \ref{sec:O4-constr} instead of the $\mathcal{O}(2)$-specific Gibbons-Hawking Ansatz which we have used in the $A_k$ case. The computation reduces effectively to the evaluation of a number of contour integrals, which we do explicitly and in a unified manner. As a check, we will then specialize the resulting ALF formulas to the case $k=0$ to retrieve the Atiyah-Hitchin metric on the moduli space of centered \mbox{charge-2} monopoles on $\mathbb{R}^3$, which plays in this case a role analogous to the one played by the Taub-NUT metric in the $A_k$ case. Finally, we compute the period matrix. 

\subsection{} 

When $\Gamma$ is the binary dihedral group of order $4k-8$, the Kleinian quotient $\mathbb{C}^2/\Gamma$ can be embedded in $\mathbb{C}^3$ as the complex surface of equation
\begin{equation}
y^2 + xz^2 = x^{k-1}
\end{equation}
with a rational double point at the origin. Its universal deformation may be described in terms of an unconstrained set of $k$ complex deformation parameters $t_l$ by the equation (see \textit{e.g.} \cite{MR1158626,Lindstrom:1999pz})
\begin{equation} \label{Dk-eq-res}
\qquad y^2 + xz^2 = \frac{1}{x}\bigg(\prod_{l=1}^k(x+t_l^2) - \prod_{l=1}^kt_l^2 \bigg) + 2z \prod_{l=1}^k t_l \rlap{.}
\end{equation}
As before, Kronheimer's classic result \cite{MR992334} assures us that we can define on the corresponding desingularized complex surface a hyperk\"ahler metric of ALE type. However, Kronheimer's approach, based on the hyperk\"ahler quotient construction, is not particularly effective when it comes to deriving explicit expressions for the metric. For this we will take, as in the $A_k$ case, the twistor road. In this case as in the $A_k$ one the basic premise of the twistor approach is that this same algebraic equation describes the complex surface not in just one but in almost all hyperk\"ahler complex structures. Accordingly, one must view the variables $x,y,z$ as well as the deformation parameters $t_l$ (which in what follows we will redenote by $t_l = i \hp x_l$) as \textit{sections} of some line bundles over the twistor space. Concretely, one assumes that  $x = x_{NS}(\zeta)$ and $x_l = x_{l,NS}(\zeta)$ are the tropical parts of real (with respect to antipodal conjugation) holomorphic global sections of the bundles $\pi^*\mathcal{O}(4)$ and $\pi^*\mathcal{O}(2)$, respectively. As for the remaining variables, their structure will be determined following a more elaborate argument which we now present. 

Notice that the deformed equation can be equivalently written as
\begin{equation}
xy^2 + \bigg(xz - \prod_{l=1}^kt_l \bigg)^{\!2} = \prod_{l=1}^k(x+t_l^2) \rlap{.}
\end{equation}
Hence, if instead of the variables $y$ and $z$ we introduce the new variables $\xi_N$ and $\xi_S$ by
\begin{equation}
\xi_{N,S} = x \hp z - \prod_{l=1}^k t_l \pm i \hp y \hp \sqrt{x} \rlap{,}
\end{equation}
then in terms of these, and with the new notation for the parameters $t_l$, the equation reads
\begin{equation}
\xi_N\xi_S = \prod_{l=1}^k (x - x_l^2) \rlap{.}
\end{equation}
In other words, what we see is that \textit{over the field of elliptic functions associated to the $\mathcal{O}(4)$ spectral curve} \eqref{O4-sc-Maj}, the equation naturally factorizes in this form. Further exchanging the new pair of variables with yet another pair, $p_N$ and $p_S$, through the relations
\begin{equation} \label{pNS-def-Dk}
\begin{aligned}
\xi_N & = e^{- 2i \hp p_N \sqrt{x} \hp \zeta} \\
\xi_S \pt & = e^{+ 2i \hp p_S \sqrt{x}/ \zeta}
\end{aligned}
\end{equation}
casts the equation in the form 
\begin{equation} \label{p-x-eq-Dk}
e^{-2i \hp (p_N\zeta - p_S/\zeta)\sqrt{x}} = \prod_{l=1}^k (x - x_l^2) \rlap{.}
\end{equation}
For our purposes it will be convenient to rewrite this formula in the following equivalent but seemingly more convoluted way:
\begin{equation} \label{pN-pS-j=2}
p_N\zeta  - \frac{p_S}{\zeta} = i \hp \frac{\partial H(x)}{\partial x} 
\end{equation}
where, by definition,
\begin{equation} \label{dHdx-Dk}
\frac{\partial H(x)}{\partial x} = \frac{1}{2\sqrt{x}} \bigg( 2\pi i + \sum_{l=1}^k \ln (x-x_l^2) \bigg) \rlap{.}
\end{equation}
The $2\pi i$ shift marks symbolically a possible logarithmic branch ambiguity, the extent of the which will be fixed subsequently based on reality and other constraints. Up to terms constant in $x$, the $H$-function is then given by 
\begin{equation} \label{H-func-Dk}
H(x) = 2\pi \hp i \hp \sqrt{x} + \sum_{l=1}^k \sum_{\pm} [(\sqrt{x} \pm x_l) \ln (\sqrt{x} \pm x_l) - (\sqrt{x} \pm x_l)] \rlap{.}
\end{equation}

The deformed quotient singularity admits a natural nowhere-vanishing holomorphic symplectic 2-form given in the $x,y,z$ variables by
\begin{equation}
\varpi = \frac{dy \wedge dz}{f_x} = \frac{dz \wedge dx}{f_y} = \frac{dx \wedge dy}{f_z}
\end{equation}
where $f(x,y,z)$ is the function whose vanishing defines the surface. Each formula is valid on a corresponding patch defined by the complement of the vanishing locus of its denominator. In the $x,\xi_N,\xi_S$ variables it may be written as
\begin{equation} \label{om-xix-Dk}
\varpi = i \pt d(\ln \xi_N) \wedge d(\sqrt{x}) = - i \pt d(\ln \xi_S) \wedge d(\sqrt{x}) 
\end{equation}
and then in the $x,p_N,p_S$ variables it assumes the form
\begin{equation}
\varpi = \frac{dp_N \wedge dx_N}{\zeta} = \zeta \pt dp_S \wedge dx_S
\end{equation}
where we have replaced the tropical form of $x$ with its arctic respectively antarctic versions.

We can now move to make the junction with the generalized Legendre transform construction. This formula makes it clear that we should identify the nowhere-vanishing holomorphic 2-form with the tropical component of the fiberwise-supported holomorphic 2-form from twistor theory, that is, 
\begin{equation}
\varpi = \omega_{NS}(\zeta)
\end{equation}
and also the $x_N,x_S$ and $p_N,p_S$ variables with the similarly denoted Legendre transform ones. The equation \eqref{pN-pS-j=2} will then correspond to the four-dimensional $j=2$ version of the gluing formula \eqref{pN-pS}, and by the same token the function \eqref{H-func-Dk} will give us precisely the generalized Legendre $H$-function characterizing ALE gravitational instantons of type~$D_k$.\footnote{\,The form of the $H$-function in the $D_k$ case was initially conjectured in \cite{MR1447294} based on expected asymptotics and the corresponding $A_k$ form, and then proven in \cite{MR1693628,MR1700937}. } These identifications provide us with an \textit{a posteriori} justification of the definitions and assumptions we have made along the way. 

\subsection{} 

The $H$-function characterizing dihedral ALF gravitational instantons, it was argued in \cite{Chalmers:1998pu, MR1693628, MR1700937}, can be obtained from its ALE counterpart \eqref{H-func-Dk} through the addition of a linear term, namely
\begin{equation} \label{H-Dk-ALF}
H^{ALF}(x) = - \alpha \hp x + H^{ALE}(x)
\end{equation}
for some real constant $\alpha$.\footnote{\,Other alternative descriptions of ALF metrics of $D_k$ type in the literature include the approach in \cite{MR2855540}, where they are constructed by means of a desingularization of the quotient of the Taub-NUT metric with the binary dihedral group of order $4k-8$ using corresponding ALE spaces, or the analytic, Monge-Amp\`ere-based approach in \cite{Auvray}.} In this case one must then consider the \textit{deformed} equation
\begin{equation} \label{xiNS-x-ALF}
\xi_N \xi_S = e^{-2\alpha\sqrt{x}} \prod_{l=1}^k (x - x_l^2) \rlap{.}
\end{equation}
Indeed, notice first that the formulas \eqref{om-xix-Dk} still give with no modification whatsoever a nowhere-vanishing holomorphic 2-form on the complex surface described by this deformed equation. Then, using the same definitions \eqref{pNS-def-Dk} for $p_N$ and $p_S$, we obtain the same equation \eqref{pN-pS-j=2} but now with the modified $H$-function \eqref{H-Dk-ALF}. The rest of the argument concerning the connection with the generalized Legendre transform construction carries over to the ALF case without any essential changes. In most of what follows we will consider a generic $\alpha$ parameter, with $\alpha=0$ corresponding to the ALE case. 

\subsection{} 

In the generalized Legendre transform approach the locally-defined twistor coordinates $p_N$ and $p_S$ are required to satisfy two important conditions:~1.~they must be well-defined at $\zeta=0$ respectively $\infty$, and 2.~they must be interchanged by antipodal conjugation up to a sign, in accordance with the second rule \eqref{real-xp}. Correspondingly, the $\xi_N$ and $\xi_S$ variables should also be well-defined at $\zeta=0$ respectively $\infty$, and should be (simply) interchanged by the action of antipodal conjugation.\footnote{\pt Note that since under antipodal conjugation $\sqrt{x_N}$ and $\sqrt{x_S}$ are interchanged up to a sign, $\sqrt{x_{NS}}$ is self-conjugate.} For the previous argument to make sense, we have to show that such solutions of the equation \eqref{xiNS-x-ALF} exist. In fact, the holomorphic structure of $x$ and $x_l$ and these two requirements determine the solutions almost entirely. To find them, we take our cue from the $A_k$ case and seek an appropriate factorization of $x_l^2-x$ and, in the ALF case, splitting of $\sqrt{x}$. 

\subsubsection{}

Let us begin with the first. Our assumptions about $x$ and $x_l$ imply that $x - x_l^2$ is the tropical part of a real (with respect to antipodal conjugation) section of $\pi^*\mathcal{O}(4)$, and so by an argument we made in the last part of section \ref{sec:Majorana-pols} it admits a factorization
\begin{equation}
x - x_l^2 = \rho_l \prod_{a=a_{l1},a_{l2}} \frac{(\zeta - a)(1+\bar{a}\zeta)}{(1+|a|^2)\zeta}
\end{equation}
where we label the roots in such a way as to have $\rho_l > 0$. Despite the obvious analogy with the $A_k$ case, this, however, turns out \textit{not} to be the factorization we need. And the reason is that $\xi_N$ and $\xi_S$\,---\,which we would like to determine by splitting up, after multiplying over all values of $l$, the resulting factors into two groups\,---\,must be elements of the field of elliptic functions associated to the $\mathcal{O}(4)$ spectral curve. As suggested in \cite{MR2177322, MR1693628, MR1700937}, a more adequate $D_k$ analogue of the $A_k$ factorization would be, rather, a factorization in terms of quasi-elliptic functions for this curve.

Guided by the meromorphic structure in $\zeta$, we consider then the alternative factorization 
\begin{equation} \label{ell-fact-xprod}
\cramped{
\prod_{l=1}^k (x - x_l^2) = \varrho \prod_a 
\frac{
\sigma(u_{\zeta} \hspace{-1pt} - u_{a}) \hp
\sigma(u_{\zeta} \hspace{-1pt} + u_{a}) \hp
\sigma(u_{\zeta} \hspace{-1pt} - u_{a^c}) \hp
\sigma(u_{\zeta} \hspace{-1pt} + u_{a^c})
}
{
\sigma(u_{\zeta} - u_0)
\sigma(u_{\zeta} + u_0)
\sigma(u_{\zeta} - u_{\infty})
\sigma(u_{\zeta} + u_{\infty}) 
} 
}
\end{equation}
in terms of Weierstrass sigma-functions, where the product is taken over all $a \in \{a_{l1},a_{l2} \, | \, l=1,\dots,k \}$ (note that these constitute only \textit{half} of all roots). The overall factor can be determined by using the ``addition theorem" \eqref{W-sgm-addthm} for the Weierstrass sigma-functions. Defining $\rho_{a_{l1}} = \rho_{a_{l2}} = \sqrt{\rho_l}$ for all $l$, we obtain
{\allowdisplaybreaks
\begin{align} \label{varrho-exprs}
\varrho 
& = \prod_a \rho_a \frac{\sigma(u_0 - u_{\infty}) \sigma(u_0 + u_{\infty})}{\sigma(u_{a} \hspace{-1pt} - u_{a^c}) \hp \sigma(u_{a} \hspace{-1pt} + u_{a^c})} \\
& = (-)^{-\frac{1}{2} \sum_a(n^{\phantom{'}}_a+n'_a)} e^{\sum_a (n^{\phantom{'}}_0\eta^{\phantom{'}}_W\upsilon'_0 + n'_0\eta'_W\upsilon^{\phantom{'}}_0 - n^{\phantom{'}}_a\eta^{\phantom{'}}_W\upsilon'_a - n'_a\eta'_W\upsilon^{\phantom{'}}_a)} |\varrho|  \rlap{.} \nonumber
\end{align}
}%
The second line displays the result as a sign times a complex phase times a positive real factor, with the last one given by
\begin{equation}
|\varrho| = \prod_a \rho_a \frac{\sigma_1(\upsilon'_0) \hp \sigma_3(\upsilon_0)}{\sigma_1(\upsilon^{\hp \prime}_{a}) \hp \sigma_3(\upsilon_{a})} e^{\frac{1}{2}(n_0^2 - n_a^2) \pt \eta^{\phantom{'}}_W\omega + \frac{1}{2}(n_0^{\prime 2} - n_a^{\prime 2}) \pt \eta'_W\omega'} > 0 \rlap{.}
\end{equation}
We employ notations established in section \ref{sec:O4-curve} and the Appendix. The Weierstrass eta-func\-tions are marked here with an index $W$ to avoid a notational clash. 

Doing the same for the Weierstrass zeta-function and defining in terms of it for any $\varphi \in \mathbb{C} \cup \{\infty\}$ the function
\begin{equation}
Z_{\varphi}(u_{\zeta}) = \frac{1}{2} [ \hp \zeta_W(u_{\zeta} \hspace{-1pt} - \hspace{-1pt} u_{\varphi}) \hspace{-1pt} + \hspace{-1pt} \zeta_W(u_{\zeta} \hspace{-1pt} + \hspace{-1pt} u_{\varphi})]
\end{equation}
then the elliptic equivalent of the splitting formula \eqref{x-split-Ak} from the $A_k$ case is 
\begin{equation} \label{ell-split-exp}
\sqrt{x} = Z_0(u_{\zeta}) - Z_{\infty}(u_{\zeta}) \rlap{.}
\end{equation}
This can be proven by using in succession the addition formula \eqref{W-zta-addthm} for the Weierstrass zeta-function and then the relation \eqref{ambig_sign_fix}. And indeed, we do have the requisite properties: $Z_0(u_{\zeta})$ is well-defined at $\zeta=\infty$, although singular at $\zeta=0$, and the other way around for $Z_{\infty}(u_{\zeta})$. Moreover, under antipodal conjugation the two component terms are not simply interchanged, rather, we have
\begin{equation} \label{Z8-ac-Z0}
\overline{Z_{\infty}(u_{\zeta^c})} = - Z_0(u_{\zeta}) + n^{\phantom{'}}_{\zeta}\eta^{\phantom{'}}_W + n'_{\zeta}\eta'_W \rlap{.}
\end{equation}
This follows from the reality property \eqref{bar-u-ac} and the monodromy property of the Weierstrass zeta-function. We will see that, far from being an inconvenience, the extra shift will turn out to play an essential role.

\subsubsection{}

Having established the two formulas \eqref{ell-fact-xprod} and \eqref{ell-split-exp}, we are now in a position to try to determine $\xi_N$ and $\xi_S$. Let us assume for $\xi_N$ the Ansatz
\begin{equation}
\xi_N  = e^{\hp 2\alpha \hp Z_{\infty}(u_{\zeta}) + 2\beta \hp u_{\zeta} + \gamma } \sqrt{|\varrho|} \prod_{a} \frac{\sigma(u_{\zeta} - u_{a}) \pt \sigma(u_{\zeta} + u_{a^c}) }{\sigma(u_{\zeta} - u_{\infty}) \sigma(u_{\zeta} + u_{\infty})} 
\end{equation}
with $\beta$ and $\gamma$ two $\zeta$-independent complex parameters to be determined. This is clearly well-defined at $\zeta=0$. The monodromy properties \eqref{W-sgm-mono}--\eqref{W-zta-mono} of the Weierstrass sigma and zeta-functions imply that under elliptic lattice shifts $u_{\zeta} \mapsto u_{\zeta} + 2 \hp \omega_i$ we have
\begin{equation}
\xi_N \mapsto \xi_N \hp e^{4 \hp \eta_{Wi}\alpha + 4 \hp \omega_i\beta - 2 \hp \eta_{Wi}\sum_a(u_{\mathstrut a}-u_{a^c})} \rlap{.}
\end{equation}
We need $\xi_N$ to be a proper doubly-periodic elliptic function, meaning that the extra factor must be trivial, which in turn only happens when the exponent is an integer multiple of $2\pi i$. Thus, requiring that at shifts by one period in the real and imaginary lattice directions the corresponding exponents be equal to $2\pi i \hp n$ respectively $2\pi i \hp n'$ for two integers $n$ and $n'$ imposes two conditions. Using the third Legendre relation \eqref{ell-Leg-rels}, these can be shown to be equivalent to
{\allowdisplaybreaks
\begin{align}
\alpha & =  n \pt \omega' - n'\omega + \frac{1}{2}\sum_a (u_a \hspace{-1pt} - u_{a^c}) \label{ell-constr-alpha} \\[-5pt]
\beta & = n'\eta^{\phantom{'}}_W - n \hp \eta'_W \rlap{.}
\end{align}
}%
By means of the antipodal reality condition on the torus the first formula can be rewritten in the form $\alpha = \sum_a \Re u_a - \big(n' + \frac{1}{2}\sum_a n^{\phantom{'}}_a \big)\omega + \big(n + \frac{1}{2}\sum_a n'_a \big) \omega'$.
Since $\omega$ is real and $\omega'$ purely imaginary, reality of $\alpha$ fixes then the value of the integer $n$ to
\begin{equation}
n = - \frac{1}{2} \sum_a n'_a 
\end{equation}
leaving us with 
\begin{equation} \label{GLT-ell-constr-0}
\sum_a \Re u_a  = \alpha + \Big(n' + \frac{1}{2}\sum_a n^{\phantom{'}}_a \Big)\omega \rlap{.}
\end{equation}
In terms of the ``enhanced" $\mathfrak{u}_a$-variables this condition can be recast, alternatively, as
\begin{equation} \label{GLT-ell-constr}
\sum_a \Re \mathfrak{u}_a = \alpha + n'\omega \rlap{.}
\end{equation}

On another hand, antipodal conjugation of the $\xi_N$ Ansatz yields 
\begin{align}
\xi_S  = (-)^{k+n+n'} e^{ \hp 2\tilde{\alpha} \hp Z_0(u_{\zeta})+ 2\tilde{\beta} \hp u_{\zeta} + \tilde{\gamma} }  \sqrt{|\varrho|} \prod_{a} \frac{\sigma(u_{\zeta} - u_{a^c})\sigma(u_{\zeta} + u_{a})}{\sigma(u_{\zeta} - u_0) \hp \sigma(u_{\zeta} + u_0)} 
\end{align}
with
{\allowdisplaybreaks
\begin{align} \label{tilded-coeffs}
\tilde{\alpha} & = - \alpha \nonumber \\
\tilde{\beta} & = - \beta \\[-2.5pt]
\tilde{\gamma} & = \bar{\gamma} +  
\sum_a [(n^{\phantom{'}}_0 \eta^{\phantom{'}}_W + n'_0 \eta'_W)( \upsilon^{\phantom{'}}_0 + \upsilon'_0) - n^{\phantom{'}}_{a} \eta^{\phantom{'}}_W \upsilon'_{a} - n'_{a} \eta'_W \upsilon^{\phantom{'}}_{a}] \rlap{.} \nonumber 
\end{align}
}%
This follows by way of the complex conjugation and monodromy properties of the Weierstrass sigma-func\-tions, the reality conditions \eqref{bar-u-ac} and \eqref{Z8-ac-Z0}, and the above corollaries of ellipticity. Note that $\xi_S$ is well-defined at $\zeta=\infty$ and, moreover, by the first two relations \eqref{tilded-coeffs} it is automatically doubly-periodic.

Finally, by multiplying the two expressions for $\xi_N$ and $\xi_S$ and then comparing the result against the elliptic splitting and factorization formulas \eqref{ell-split-exp} and \eqref{ell-fact-xprod}, with $\varrho$ given by the second expression \eqref{varrho-exprs}, we conclude that  the equation \eqref{xiNS-x-ALF} holds if and only if
\begin{equation}
\Re \gamma = -k \hp (n^{\phantom{'}}_0\eta^{\phantom{'}}_W\upsilon^{\phantom{'}}_0 + n'_0\eta'_W\upsilon'_0)
\end{equation}
and
\begin{equation}
n' = - \frac{1}{2}\sum_an_a + k + 2\hp s 
\end{equation}
for some $s \in \mathbb{Z}$. By the very existence of a solution our choice of Ansatz is vindicated \textit{a posteriori}. The fact that only the real part of the parameter $\gamma$ is fixed is consistent with $\xi_N$ and $\xi_S$ being defined only up to opposite constant complex phases. 

The integer $s$ can be understood as a compensator or a ``sink" for two types of ambiguities: the $(\mathbb{Z}_2)^k$-ambiguity associated to the choice of the $a$ roots,\footnote{\pt For each $l = 1,\dots,k$ the condition $\rho_l >0$ fixes, as we have seen, a relative ambiguity in the choice of $a_{l1}$ and $a_{l2}$, but one is still left with a residual $\mathbb{Z}_2$ ambiguity deriving from the fact that if $a_{l1}$ and $a_{l2}$ give $\rho_l >0$ then so do their antipodal conjugates.} and the $\mathbb{Z}$-ambiguity associated to the choice of fundamental domain. On the other hand, the integer $n'$ is not affected by either of these, as can be seen for instance from the form \eqref{GLT-ell-constr} of the first ellipticity constraint in the light of the properties of the $\mathfrak{u}_a$-variables discussed in $\S$~\ref{ssec:ac-torus}.

We can decouple the $l$-th monopole by taking the moduli $\vec{r}_l$ of the corresponding $\mathcal{O}(2)$ section $x_l$ to infinity \cite{Seiberg:1996nz, Sen:1997kz, Chalmers:1998pu}. The roots of $x_N - x_{lN}^2$ approach then those of $x_{lN}$, and in particular it follows that $a_{l1}$ and $a_{l2}$ approach antipodally-conjugated values. In this limit, therefore, $n_{a_{l1}} = n_{a_{l2}}$ and $\Re u_{a_{l1}} + \Re u_{a_{l2}} = n_{a_{l1}}\omega$, and it is clear that the $l$-terms in the two sums in \eqref{GLT-ell-constr-0} cancel each other out, leaving behind precisely the constraint for the $D_{k-1}$ space. For $\alpha \neq 0$ we may continue to remove more and more monopoles until we are left with none, which will leave the constraint equation for the Atiyah-Hitchin space (see $\S$~\ref{ssec:AH} ahead). This prompts us to take in this case $n'=2$. 

\subsection{} 

Next, we pursue the explicit computation of the hyperk\"ahler $D_k$ metrics and symplectic forms. The input data for the calculation consists of the holomorphic $H$-function, which in the ALE case is given by the formula \eqref{H-func-Dk} and in the ALF one by the formula \eqref{H-Dk-ALF}, together with a choice of contour for the real-valued contour integral $L$. The relevant generalized Legendre transform framework is the $\mathcal{O}(4)$ or $j=2$ one, and so the apposite construction is the quaternionic Ansatz of Proposition \ref{O4-Ansatz}. In this, the metric and symplectic forms are effectively determined by the second derivatives $L_{x_0x_m}$ with $m=-2,\dots,2$. In addition, one also has to impose the characteristic differential constraint $L_{x_0} = 0$. The computations quickly reduce to evaluating two basic kinds of elliptic integrals, and the approach best suited to this Ansatz is to try to express them in terms of Weierstrass-type elliptic and quasi-elliptic functions. This is because we want if not the transcendental functions themselves then at least their coefficients to depend as explicitly as possible on the Majorana coefficients, and through them, on the quaternionic variables.

Let us begin by examining the constraint. From the formula \eqref{dHdx-Dk}, taking also into account $\alpha$-deformations, we obtain
{\allowdisplaybreaks
\begin{align}
L_{x_0} & = \oint \frac{d\zeta}{2\pi i \zeta} \frac{\partial H(x)}{\partial x} \\
& = - \alpha \oint_{2\hp C_0} \frac{d\zeta}{2\pi i \zeta} + \oint_{n\hp C'-n'C} \frac{d\zeta}{2\sqrt{x_N(\zeta)}} + \sum_a \int_{\hspace{-0.5pt}a^c}^{\pt a} \! \frac{d\zeta}{2\sqrt{x_N(\zeta)}} \nonumber \\
& = - 2\hp \alpha + 2\hp n \pt \omega' - 2\hp n'\omega + \sum_a(u_a-u_{a^c}) \rlap{.} \nonumber
\end{align}
}%
The contour consists of separate components for each term. For the $\alpha$-dependent term we consider a simple contour wrapped counter-clockwise around 0 twice\,---\,or, equivalently, a contour wrapping once around both $0$ and $\infty$ with opposite orientations. For the second term we take a contour which wraps $n$ times counter-clockwise and $n'$ times clockwise around the elliptic cycles corresponding to the real and imaginary periods of the $\mathcal{O}(4)$ spectral curve, respectively. For the logarithmic terms in the sum we choose figure eight-shaped contours circling pairs of antipodally-conjugated roots $a$ and $a^c$. The last line follows after a further change to complex torus variables as in equation \eqref{1-f-MWJ}. So in the end what we find, remarkably, is that the generalized Legendre transform constraint $L_{x_0} = 0$ is precisely the same as the ellipticity constraint \eqref{ell-constr-alpha}! \cite{MR1693628, MR1700937}

To compute $L_{x_0x_m}$, we differentiate again the second line above with respect to $x_m$. The $\alpha$-term drops out completely, and the terms under the sum give 
\begin{align} \label{diff-int}
\frac{\partial}{\partial x_m} \int_{\hspace{-0.5pt}a^c}^{\pt a} \!\frac{d\zeta}{2\sqrt{x_N(\zeta)}} = & \
\frac{1}{2\sqrt{x_N(a)}} \frac{\partial a}{\partial x_m} - \frac{1}{2\sqrt{x_N(a^c)}} \frac{\partial a^c}{\partial x_m} \\
& - \frac{1}{2} \int_{\hspace{-0.5pt}a^c}^{\pt a} \! \frac{d\zeta}{2\sqrt{x_N(\zeta)}} \frac{\zeta^{2-m}}{x_N(\zeta)} \rlap{.} \nonumber
\end{align}
The remaining term leads to a similar but simpler variant of this formula, with a closed integration contour and hence no boundary terms. For the evaluation of the resulting integrals in terms of Weierstrass elliptic and quasi-elliptic functions we rely on the following key formula:
\begin{equation} \label{indef-int}
\int \frac{d\zeta}{2\sqrt{x_N(\zeta)}} \frac{\zeta^{2-m}}{x_N(\zeta)}
= (-)^m A_{-m} u_{\zeta} - (-)^m B_{-m} Z_{\infty}(u_{\zeta}) + \frac{L^{(x_N)}_{1-m}(\zeta)}{\sqrt{x_N(\zeta)}} \ + \, \mathscr{C} \rlap{.}
\end{equation}
Before we consider its proof, we will first take some time to define the coefficients $A_{-m}$, $B_{-m}$ and the functions $L^{(x_N)}_{1-m}(\zeta)$, and discuss their properties. Let us just mention at this point that $A_{-m}$ and $B_{-m}$ denote certain homogeneous polynomials in $\mathbb{Q}[x_{-2},\dots,x_{+2}]$ of degree five and four, respectively, and $L^{(x_N)}_{1-m}(\zeta)$ stands for a polynomial of order three in $\zeta$ with coefficients given by degree-five homogeneous polynomials in $\mathbb{Q}[x_{-2},\dots,x_{+2}]$. 

\subsubsection{} 

Explicitly, the coefficients of the quasi-elliptic functions are defined by 
{\allowdisplaybreaks
\begin{equation}
\begin{aligned}
A_m & = (-)^m \frac{1}{\Delta} \bigg( 2g_2^2 \frac{\partial g_2}{\partial x_{-\rlap{$\scriptstyle m$}}} \ - 9g_3 \frac{\partial g_3}{\partial x_{-\rlap{$\scriptstyle m$}}} \hspace{6pt} \bigg) = \phantom{+} (-)^m \frac{1}{6} \frac{\partial \ln \Delta}{\partial x_{-m}} \\
B_m & = (-)^m \frac{1}{\Delta} \bigg( 9g_3 \frac{\partial g_2}{\partial x_{-\rlap{$\scriptstyle m$}}} \ - 6g_2 \frac{\partial g_3}{\partial x_{-\rlap{$\scriptstyle m$}}} \hspace{6pt} \bigg) \hp = - (-)^m \frac{g_2}{9g_3} \frac{\partial \ln J}{\partial x_{-m}} \\
\end{aligned}
\end{equation}
}%
where $J = 4g_2^3/(4g_2^3 - 27g_3^2)$ is Klein's absolute invariant. In view of the formulas \mbox{\eqref{r2r3}--\eqref{Weierstrass-coeffs}} for $g_2,g_3$ they are clearly expressible as homogeneous polynomials in the Majorana coefficients. Moreover, and most relevantly for us, they satisfy alternating reality conditions and have spherical tensor properties. Specifically, at an $SO(3)$ transformation of the type \eqref{spinor-rot} with $j=2$ they transform covariantly according to the rules
{\allowdisplaybreaks
\begin{equation}
\begin{aligned}
A^2_m & = \sum_{m'=-2}^2 D^2_{mm'}(q) \pt \mathring{A}^2_{m'} \\[-3pt]
B^2_m & = \sum_{m'=-2}^2 D^2_{mm'}(q) \pt \mathring{B}^2_{m'}
\end{aligned}
\end{equation}
}%
where we normalize by $A^2_m = (c^2_{m})^{-1} A^{\phantom{2}}_m$ and $B^2_m = (c^2_{m})^{-1} B^{\phantom{2}}_m$. The same principle is at work here which underlies the discussion surrounding equation~\eqref{dgdx}. In particular, if we take $\mathring{x}_m$ to be of the form \eqref{ref-sp-j=2}, we then get 
\begin{equation} \label{AB-ring}
(\hp \mathring{A}_m \hp ) = - \frac{1}{\sqrt{\Delta}}
\begin{pmatrix}
e_1^2+4e_1e_3+e_3^2 \\
0 \\
e_1^2-e_3^2 \\
0 \\
e_1^2+4e_1e_3+e_3^2
\end{pmatrix}
\quad \text{and} \quad
(\hp \mathring{B}_m \hp ) =  \frac{1}{\sqrt{\Delta}}
\begin{pmatrix}
3(e_1+e_3) \\
0 \\
- (e_1-e_3) \\
0 \\
3(e_1+e_3)
\end{pmatrix}
\rlap{.}
\end{equation}
Other properties of some practical use include the relations
\begin{equation}
\begin{aligned}
& 4B_{m-2}x_{+2} - 3B_{m-1}x_{+1} + 2B_mx_0 - B_{m+1}x_{-1} = 0 && (m=0,1) \\
& 6B_{m-2}x_{+2} - 3B_{m-1}x_{+1} + \phantom{2}B_mx_0  = - 3 A_m && (m=0,1,2)
\end{aligned}
\end{equation}
and their complex conjugates.

\subsubsection{} 

To any quartic Majorana polynomial $x_N(\zeta)$, regarding its roots $a_1,\dots,a_4$ as functions of the Majorana coefficients, we associate the cubic polynomials
\begin{equation}
L^{(x_N)}_{1-m}(\zeta) = \sum_{i=1}^4 \ell_i(\zeta) \frac{\partial a_i}{\partial x_m}
\end{equation}
with $m \in \{-2,\dots,2\}$, defined by means of the Lagrange interpolation polynomials
\begin{equation}
\ell_i(\zeta) = \prod_{\shortstack{$\scriptstyle j=1$ \\ $\scriptstyle j \neq i$}}^4 \frac{\zeta-a_j}{a_i-a_j} \rlap{.}
\end{equation}
For any such integer $m$, these are the unique polynomials of order at most three interpolating through the first derivatives of the roots with respect to $x_m$. Indeed, we have $\ell_i(a_j) = \delta_{ij}$, and so at the roots 
\begin{equation} \label{L@roots}
L^{(x_N)}_{1-m}(a_i) =  \frac{\partial a_i}{\partial x_m} \rlap{.}
\end{equation}
Implicit differentiation gives us for the derivatives of the roots the expressions
\begin{equation}
\frac{\partial a_i}{\partial x_m} = - \frac{1}{x_{-2}}\frac{a_i^{2-m}}{\prod_{j=1, j\neq i}^4 (a_i-a_j)} \rlap{.}
\end{equation}
In view of this formula, the domain of $m$ can be analytically continued to include integer values for which $x_m$ is not defined. In particular, for $m \in \{-1,\dots,4\}$ we have the following properties:
{\allowdisplaybreaks
\begin{equation} \label{Lag-interp-prop}
\begin{aligned}
1.\ \ & L^{(x_N)}_m(\zeta) = \zeta L^{(x_N)}_{m-1}(\zeta) + (-)^{m-2}B_{m-2} \hp x_N(\zeta) \\[0pt]
2.\ \ & L^{(x_N)}_m(\zeta) = \sum_{i=0}^3 \Big( \sum_{\shortstack{$\scriptstyle r+s+i = m$ \\ $\scriptstyle r+s \geq 0 \hfill$ \\$\scriptstyle \phantom{r+}s+i \geq 2 \pt\pt\hp$}} \hspace{-10pt} (-)^r B_{r}x_{s} \hspace{3pt} - \hspace{-8pt} \sum_{\shortstack{$\scriptstyle r+s+i = m$ \\ $\scriptstyle r+s < 0 \hfill$ \\$\scriptstyle \phantom{r+}s+i < 2 \pt\pt\hp$}} \hspace{-10pt} (-)^r B_{r}x_{s} \Big) \zeta^i  \\
3.\ \ & \overline{L^{(x_N)}_m(\zeta^c)} = \frac{(-)^{\rlap{$\scriptstyle 3-m$}}}{\zeta^3} \hspace{14pt} L^{(x_N)}_{3-m}(\zeta) \rlap{.}
\end{aligned}
\end{equation}
}%
The first one is a recurrence relation. The second one states that, remarkably, the polynomial coefficients may be expressed entirely in terms of the Majorana coefficients. Finally, the third property is a reality condition reflecting the one satisfied by the Majorana polynomial. \\


With these definitions in place, we are now able to return to the indefinite integral \eqref{indef-int} to discuss its proof. But rather than taking the reader through the steps by which this formula was obtained, we simply check here that the derivative of the result yields back the integrand. The crucial facts to notice for this are the following:
\begin{equation}
\begin{aligned}
\frac{\partial \hp u_{\zeta}}{\partial \zeta} & = \frac{1}{2\sqrt{x_N(\zeta)}} \\[-2pt]
\frac{\partial Z_{\infty}(u_{\zeta})}{\partial \zeta} & = - \frac{1}{12\sqrt{x_N(\zeta)}} \frac{\partial^2 x_N(\zeta)}{\partial \zeta^2} \rlap{.}
\end{aligned}
\end{equation}
The first derivative is just a rephrasing of the differential relation \eqref{1-f-MWJ}. To obtain the second one we first used the addition formula \eqref{W-zta-addthm} and then the relation \eqref{ambig_sign_fix} to rewrite $Z_{\infty}(u_{\zeta})$. Derivating the resulting expression and making use of the first derivative, one can eventually cast the outcome in the form above. With these two differentiation formulas, checking the integral \eqref{indef-int} reduces then to a straightforward algebraic computation.

Having established the integration formula \eqref{indef-int}, we can now deploy it to evaluate $L_{x_0x_m}$, which is, recall, our original objective. A swift calculation gives us
{\allowdisplaybreaks
\begin{align}
L_{x_0x_m} = & - \frac{1}{2} (-)^mA_{-m} \big[\sum_a(u_a - u_{a^c}) + 2\hp n\hp \omega' - 2\hp n'\omega \hp\big] \\[-3pt]
&  + \frac{1}{2} (-)^mB_{-m} \big[\sum_a(Z_{\infty}(u_a) - Z_{\infty}(u_{a^c})) + 2\hp n\hp \eta'_W - 2\hp n'\eta^{\phantom{'}}_W \big] \nonumber \\
& + \sum_a [C_m(a)-C_m(a^c)] \nonumber
\end{align}
}%
where, by definition,
\begin{equation} \label{alg-term-fct}
C_m(\zeta) = \frac{L^{(x^{\phantom{2}}_N-x_{lN}^2)}_{1-m}(\zeta) - L^{(x_N)}_{1-m}(\zeta)}{2\sqrt{x_N(\zeta)}} \rlap{.}
\end{equation}
The cubic $L$-polynomials associated to the quartic Majorana polynomial $x^{\phantom{2}}_N-x_{lN}^2$ (for which, remember, all $a$ and $a^c$ are roots) arise by way of the corresponding property \eqref{L@roots} from the boundary terms in \eqref{diff-int}. 

For any root $a$ of $x^{\phantom{2}}_N-x_{lN}^2$, based on the third and then the first property \eqref{Lag-interp-prop}, we have 
\begin{equation}
\overline{C_m(a^c)} = -(-)^m C_{-m}(a) + \frac{1}{2} B_m\sqrt{x_{NS}(a)} \rlap{.}
\end{equation}
Note that the shift comes entirely from the second term in \eqref{alg-term-fct}. 

In the same way we have introduced an ``enhanced" version $\mathfrak{u}_{\zeta}$ of $u_{\zeta}$, let us introduce also an ``enhanced" version of $Z_{\infty}(u_{\zeta})$ by $\mathfrak{Z}(u_{\zeta}) = Z_{\infty}(u_{\zeta}) - (n_{\zeta}/2)\hp\eta^{\phantom{'}}_W - (n'_{\zeta}/2)\hp\eta'_W$. From the properties \eqref{ell-split-exp} and \eqref{Z8-ac-Z0} it follows that this satisfies the reality condition
\begin{equation}
\overline{ \mathfrak{Z}(u_{\zeta^c}) } = - \mathfrak{Z}(u_{\zeta}) - \sqrt{x(\zeta)}
\end{equation}
valid for any $\zeta \in \mathbb{C}^{\times}$.

By using these two reality properties and new definition the expression above can be recast in the form
{\allowdisplaybreaks
\begin{align}
L_{x_0x_m} = & - (-)^mA_{-m} \big[\sum_a \Re \mathfrak{u}_a - n' \omega \hp\big] \\[-2pt]
& + (-)^mB_{-m} \big[\sum_a \Re \mathfrak{Z}(u_a) - n' \eta^{\phantom{'}}_W \big] \nonumber \\
& + \sum_a [C_m(a) + (-)^m\overline{C_{-m}(a)} \hp] \rlap{.} \nonumber
\end{align}
}%
Taking furthermore into account the constraint \eqref{GLT-ell-constr} and recalling that $L_{x_0x_m} =\bar{h}_m$, we arrive in the end at the following result: 

\begin{proposition} \label{Dk-metric}
The Ansatz of Proposition \ref{O4-Ansatz} with
\begin{equation} \label{hm-Dk}
h_m = - A_m \hp \alpha + B_m \big[\sum_a \emph{\Re} \mathfrak{Z}(u_a) - n' \eta^{\phantom{'}}_W \big] + \sum_a [(-)^mC_{-m}(a) + \overline{C_m(a)} \hp]
\end{equation}
(\pt$m=-2,\dots,2$) and the constraint
\begin{equation}
\sum_a \emph{\Re} \mathfrak{u}_a = \alpha + n'\omega
\end{equation}
describes a hyperk\"ahler metric and corresponding symplectic forms on an ALE ($\alpha=0$) or ALF ($\alpha \neq 0$) gravitational instanton of type $D_k$. 
\end{proposition}

\subsection{The Atiyah-Hitchin or ALF $D_0$ case} \label{ssec:AH}

Let us specialize now these formulas to the ALF case with $k=0$. Based on generic considerations, this should give us the metric of Atiyah and Hitchin on the moduli space of centered charge-2 monopoles \cite{MR934202}, which leads the series for ALF gravitational instantons of type $D_k$ similarly to how the Taub-NUT metric does for the ALF gravitational instantons of type $A_k$. Taking $n'=2$, the constraint reads
\begin{equation}
\alpha = - 2\hp \omega
\end{equation}
while the expression \eqref{hm-Dk} becomes simply
\begin{equation}
h_m = - \alpha A_m \rlap{.}
\end{equation}
Reverting to ringed quantities and using for $\mathring{A}_m$ the form \eqref{AB-ring} we obtain immediately the expressions
\begin{equation}
r_2(\mathring{h}) = \frac{4\hp\alpha^2g_2^2}{\Delta} 
\qquad\quad
r_3(\mathring{h}) = \frac{4\hp\alpha^3g_3}{\Delta} 
\end{equation}
as well as a set of explicit formulas for all $\mathring{H}_m$. In particular, we have
\begin{equation}
\mathring{H}'_0= -  \frac{\alpha^2 \rho^2e_2^2}{\Delta}
\end{equation}
which is, as required, non-positive. Substituting into the Ansatz yields after straightforward algebraic manipulations the following formula for the metric
\begin{equation}
g = \frac{e_1e_3}{8\hp e_2} \frac{d\rho^2}{\rho^2} + \frac{2\hp e_2e_1}{e_3}\sigma_1^2 + \frac{2\hp e_3e_2}{e_1}\sigma_2^2 + \frac{2\hp e_1e_3}{e_2}\sigma_3^2
\end{equation}
and moreover
\begin{equation}
\begin{aligned}
\mathring{\omega}_1 & = - 4\hp e_3\pt \sigma_2 \wedge \sigma_3 + e_1\hp \sigma_1 \wedge \frac{d\rho}{\rho} \\[-2pt]
\mathring{\omega}_2 & = \phantom{+} 4\hp e_1\hp \sigma_3 \wedge \sigma_1 \hp - e_3\pt \sigma_2 \wedge \frac{d\rho}{\rho} \\[-2pt]
\mathring{\omega}_3 & = \phantom{+} 4\hp e_2\pt \sigma_1 \wedge \sigma_2 - \frac{e_1e_3}{e_2}\hp  \sigma_3 \wedge \frac{d\rho}{\rho} \rlap{,}
\end{aligned}
\end{equation}
from which the hyperk\"ahler symplectic forms may be simply obtained by a further rotation as described in Proposition \ref{O4-Ansatz}. We have dropped in each case from the formulas an overall factor of $-2\hp\alpha$. As one can easily ascertain, in this particular case the formulas of Proposition \ref{Dk-metric} do reproduce indeed just as expected the Atiyah-Hitchin metric. 

In the same vein, the ALF $D_1$ and $D_2$ metrics we have found should coincide, presumably, with the metrics on the double-cover of the Atiyah-Hitchin manifold and on the minimal resolution of $(\mathbb{R}^3 \times S^1)/\mathbb{Z}_2$ from \cite{MR2855540} and \cite{MR1255427} (see also \cite{MR757213}), respectively.

\subsection{The period matrix} 

In the last installment of our discussion we return to the general ALE case to consider the period matrix. Given a gravitational instanton $M$, its period matrix is obtained by pairing a basis of exceptional homology classes in $H_2(M,\mathbb{Z})$ with the cohomology classes in $H^2(M,\mathbb{R})$ determined by the hyperk\"ahler symplectic forms. The exceptional divisor of the minimal resolution of the Kleinian quotient $\mathbb{C}^2/\Gamma$ with $\Gamma$ the binary dihedral group of order $4k-8$ consists of $k$ rational curves with intersection matrix equal to minus the Cartan matrix corresponding to a Dynkin diagram of type $D_k$. Following the pattern we laid out in the $A_k$ case, let us try now to describe somewhat more explicitly a set of representative 2-cycles in their homology classes. 

One can regard the universal deformation equation \eqref{Dk-eq-res} as a fibration of quadrics over the $x$-plane, with singular fibers at the points $x=-t_l^2=x_l^2$ (where the symmetric matrix defining the quadric has vanishing determinant). In this complex fibration, for any pair of distinct indices $l_1,l_2 \in \{1,\dots,k\}$, we consider a real 2-cycle $\gamma_{\hp l_1l_2}$ in the form of a circle fibration over a 1-dimensional path in the $x$-plane connecting the points $x_{l_1}^2$ and $x_{l_2}^2$, with the circle degenerating at the endpoints. In terms of the $x$ and $p$ variables, such a 2-cycle can be parametrized as follows:
{\allowdisplaybreaks
\begin{equation}
\begin{aligned}
& \text{$x  = x(\lambda)$ such that $x(0)=x_{l_2}$ and $x(1)=x_{l_1}$} \\
& p_N = \frac{\phi}{2\sqrt{x_N(\lambda)}} + p_N^0(\lambda) \\
& p_S = \pt \frac{\phi}{2\sqrt{x_S(\lambda)}} \pt + p_S^0(\lambda)
\end{aligned}
\end{equation}
}%
where $\phi \in [0,2\pi]$ is the circle parameter and $\lambda \in [0,1]$ the path parameter. 

But, remember, the equation \eqref{Dk-eq-res} describes the graviton in a hyperk\"ahler complex structure labeled by $\zeta$, and the question is, does this parametrization of the 2-cycle break if we change $\zeta$? This is in no way a trivial question. In the $A_k$ case the three components $x_m$ of $x$ could be considered in principle independent variables and, because of that, one could make the case that the answer was ``no" easier. In the $D_k$ case, however, there are five components and they are related by a fundamental generalized Legendre transform constraint, and hence cannot be considered as independent. So \textit{a priori} it is by no means guaranteed that if in a given complex structure we choose $x$ to be parametrized only by $\lambda$, or $p_N$ and $p_S$ to be parametrized in the above simple way in terms of $\phi$, then the same thing will happen in other complex structures. The important fact to keep in mind at this point is that all the generalized Legendre transform equations, including the constraint, are fully encapsulated in the equation \eqref{p-x-eq-Dk} through the powerful means of holomorphicity in $\zeta$. Then there are two things to notice: First, the parametrization of $p_N$ and $p_S$ is compatible with antipodal conjugation. Second, the parameter $\phi$ cancels out in the exponent on the right-hand side of the equation \eqref{p-x-eq-Dk}. And once it cancels out, the remaining dependence on $\lambda$ is consistent and exclusive order by order in $\zeta$ in the remaining equation. This implies in particular that the basic structure of the parametrization does not change at a change of $\zeta$. Or, otherwise put, our choice of 2-cycle does not depend on $\zeta$. 

Using the complex Darboux expression for the twisted symplectic form we have then, successively, 
{\allowdisplaybreaks
\begin{align}
\int_{\gamma_{\hp l_1l_2}} \! \omega_N(\zeta) 
& = \int_{\gamma_{\hp l_1l_2}} dp_N(\zeta) \wedge dx_N(\zeta) \\[2pt] 
& = \int_{\gamma_{\hp l_1l_2}} d\phi \wedge \frac{dx_N(\zeta)}{2\sqrt{x_N(\zeta)}} \nonumber \\ 
& = \int_0^{2\pi} \!\! d\phi \int_{x_{l_2N}(\zeta)^2}^{\mathstrut  x_{l_1N}(\zeta)^2} \frac{dx_N(\zeta)}{2\sqrt{x_N(\zeta)}} \nonumber \\[5pt]
& = 2\pi \hp [\sqrt{x_{\raisebox{0pt}[\height][0pt]{$\scriptstyle l_1$}N}(\zeta)^2} - \sqrt{x_{\raisebox{0pt}[\height][0pt]{$\scriptstyle l_2$}N}(\zeta)^2} \hp] \nonumber
\end{align}
}%
The square roots in the last line can be trivially extracted to yield a linear superposition of $x_{l_1N}(\zeta)$ and $x_{l_2N}(\zeta)$, although with possibly ambiguous signs. Since $\zeta$ is arbitrary, we may eventually conclude that it is possible to choose a subset of $k$ elements from the set of 2-cycles of the form $\gamma_{\hp l_1l_2}$\,---\,which we rename $\gamma_i$, with $i =1,\dots,k$\,---\,such that their homology classes span $H^2(M,\mathbb{Z})$ and 
\begin{align}
\int_{\gamma_i} \vec{\omega} = 
\begin{cases}
2\hp \pi \hp (\vec{r}_i \hp - \hp \vec{r}_{i+1})  & \text{for $i = 1,\dots,k-1$} \\[3pt]
2\hp \pi \hp (\vec{r}_{k-1} \! + \vec{r}_k) & \text{for $i = k$} 
\end{cases}
\end{align}
where the $\vec{r}_l$ are the moduli of the $\mathcal{O}(2)$ sections $x_l$. The Torelli-type theorem of \cite{MR992335} asserts that the periods of the three hyperk\"ahler symplectic forms on the exceptional cycles determine the metric uniquely up to diffeomorphism. Thus, for each set of parameters $\vec{r}_i$, the metric we have described can be regarded as a representative of this unique equivalence class.

\appendix

\section{Elliptic formulas}

\setcounter{equation}{0}\renewcommand{\theequation}{\Alph{section}.\arabic{equation}}

For ease of reference, we collect in this appendix a short compendium of useful elliptic formulas. 

The so-called addition theorem for the Weierstrass sigma-function is in fact the multiplicative formula 
\begin{equation} \label{W-sgm-addthm}
\frac{\sigma(u-v)\sigma(u+v)}{\sigma(u)^2\sigma(v)^2} = \wp(v) - \wp(u) \rlap{.}
\end{equation}
Logarithmic differentiation with respect to the variable $u$ yields the following proper addition formula for the Weierstrass zeta-function:
\begin{equation} \label{W-zta-addthm}
\zeta(u-v) + \zeta(u+v) = 2 \hp \zeta(u) + \frac{\wp'(u)}{\wp(u)-\wp(v)} \rlap{.}
\end{equation}
The monodromy properties for the sigma and zeta-functions are given by
\begin{align}
\sigma(u + 2\hp \omega_i) & = - e^{2 \hp \eta_{i} (u + \omega_i)} \sigma(u) \label{W-sgm-mono} \\[1pt]
\zeta(u + 2 \hp \omega_i) & = \zeta(u) + 2 \hp \eta_{i} \rlap{,} \label{W-zta-mono}
\end{align}
with $i=1,2,3$, where $\omega_i$ satisfying $\omega_1+\omega_2+\omega_3=0$ are the Weierstrass half-periods and $\eta_i = \zeta(\omega_i)$ are the Weierstrass eta-functions. These satisfy Legendre's relations
\begin{equation} \label{ell-Leg-rels}
\eta_3\hp \omega_2 - \eta_2\hp \omega_3 = \eta_2\hp \omega_1 - \eta_1\hp \omega_2 = \eta_1\hp \omega_3 - \eta_3 \hp \omega_1 = \frac{i\pi}{2} \rlap{.}
\end{equation}
In addition, we have the differential formulas
\begin{equation}
d\omega_i = - \frac{2g_2^2dg_2 - 9g_3dg_3}{2\Delta} \omega_i + \frac{9g_3dg_2 - 6g_2dg_3}{2\Delta} \eta_{i} \rlap{.}
\end{equation}
where $g_2$, $g_3$ are the Weierstrass coefficients and $\Delta= 4 g_2^3 - 27g_3^2$ is the Weierstrass discriminant. (\textit{Caveat lector}: compared to standard references, in this paper we are using a slightly unusual normalization conventions for these.) Finally, three \textit{associated} Weierstrass sigma-functions are defined by
\begin{equation}
\sigma_i(u) = e^{-\eta_iu}\frac{\sigma(u+\omega_i)}{\sigma(\omega_i)}
\end{equation}
and satisfy, among other things, the property
\begin{equation} \label{assoc-sigma-propr}
\wp(u) - e_i = \frac{\sigma_i(u)^2}{\sigma(u)^2} \rlap{,}
\end{equation}
where the $e_i$ denote the Weierstrass roots.

In accordance with the usual conventions we sometimes denote $\omega_1$ and $\omega_3$ by $\omega$ and $\omega'$ and, also, $\eta_1$ and $\eta_3$ by $\eta$ and $\eta'$. Note that in the main text we append a $W$ subscript to the usual symbols of the Weierstrass zeta and eta-functions in order to differentiate them from other similarly denoted quantities.

\section*{Acknowledgements}

This work has benefited significantly from many discussions with Martin Ro\v{c}ek, to whom the author wishes to thank for his continued advice and support. The presentation owes also many improvements to Andrew Hanson's spirited probing questions.

\bigskip 

\Addresses

\bibliography{HK-LT}

\begin{thebibliography}{10}

\bibitem{MR934202}
M.~Atiyah and N.~Hitchin.
\newblock {\em The geometry and dynamics of magnetic monopoles}.
\newblock M. B. Porter Lectures. Princeton University Press, Princeton, NJ,
  1988.

\bibitem{Auvray}
H.~Auvray.
\newblock From {ALE} to {ALF} gravitational instantons.
\newblock 2012.
\newblock arXiv:1210.1654.

\bibitem{MR1848654}
R.~Bielawski.
\newblock Twistor quotients of hyperk\"ahler manifolds.
\newblock In {\em Quaternionic structures in mathematics and physics ({R}ome,
  1999)}, pages 7--21 (electronic). Univ. Studi Roma ``La Sapienza'', Rome,
  1999.

\bibitem{MR2855540}
O.~Biquard and V.~Minerbe.
\newblock A {K}ummer construction for gravitational instantons.
\newblock {\em Comm. Math. Phys.}, 308(3):773--794, 2011.

\bibitem{Chalmers:1998pu}
G.~Chalmers, M.~Ro\v{c}ek, and S.~Wiles.
\newblock {Degeneration of ALF {$D_n$} metrics}.
\newblock {\em JHEP}, 01:009, 1999.

\bibitem{MR2177322}
S.~A. Cherkis and N.~J. Hitchin.
\newblock Gravitational instantons of type {$D_k$}.
\newblock {\em Comm. Math. Phys.}, 260(2):299--317, 2005.

\bibitem{MR1693628}
S.~A. Cherkis and A.~Kapustin.
\newblock {$D_k$} gravitational instantons and {N}ahm equations.
\newblock {\em Adv. Theor. Math. Phys.}, 2(6):1287--1306 (1999), 1998.

\bibitem{MR1700937}
S.~A. Cherkis and A.~Kapustin.
\newblock Singular monopoles and gravitational instantons.
\newblock {\em Comm. Math. Phys.}, 203(3):713--728, 1999.

\bibitem{MR1255427}
A.~S. Dancer.
\newblock Nahm's equations and hyper-{K}\"ahler geometry.
\newblock {\em Comm. Math. Phys.}, 158(3):545--568, 1993.

\bibitem{Gibbons:1979zt}
G.~W. Gibbons and S.~W. Hawking.
\newblock {Gravitational Multi-Instantons}.
\newblock {\em Phys. Lett.}, B78:430, 1978.

\bibitem{hanson2006}
A.~Hanson.
\newblock {\em Visualizing Quaternions}.
\newblock The Morgan Kaufmann Series in Interactive 3D Technology. Elsevier
  Science, 2006.

\bibitem{MR3543182}
H.-J. Hein and C.~LeBrun.
\newblock Mass in {K}\"ahler geometry.
\newblock {\em Comm. Math. Phys.}, 347(1):183--221, 2016.

\bibitem{MR520463}
N.~J. Hitchin.
\newblock Polygons and gravitons.
\newblock {\em Math. Proc. Cambridge Philos. Soc.}, 85(3):465--476, 1979.

\bibitem{MR623721}
N.~J. Hitchin.
\newblock K\"ahlerian twistor spaces.
\newblock {\em Proc. London Math. Soc. (3)}, 43(1):133--150, 1981.

\bibitem{MR757213}
N.~J. Hitchin.
\newblock Twistor construction of {E}instein metrics.
\newblock In {\em Global {R}iemannian geometry ({D}urham, 1983)}, Ellis Horwood
  Ser. Math. Appl., pages 115--125. Horwood, Chichester, 1984.

\bibitem{MR887284}
N.~J. Hitchin.
\newblock The self-duality equations on a {R}iemann surface.
\newblock {\em Proc. London Math. Soc. (3)}, 55(1):59--126, 1987.

\bibitem{MR877637}
N.~J. Hitchin, A.~Karlhede, U.~Lindstr{\"o}m, and M.~Ro{\v{c}}ek.
\newblock Hyper-{K}\"ahler metrics and supersymmetry.
\newblock {\em Comm. Math. Phys.}, 108(4):535--589, 1987.

\bibitem{MR1447294}
I.~T. Ivanov and M.~Ro{\v{c}}ek.
\newblock Supersymmetric {$\sigma$}-models, twistors, and the
  {A}tiyah-{H}itchin metric.
\newblock {\em Comm. Math. Phys.}, 182(2):291--302, 1996.

\bibitem{MR1669956}
D.~Kaledin and M.~Verbitsky.
\newblock Non-{H}ermitian {Y}ang-{M}ills connections.
\newblock {\em Selecta Math. (N.S.)}, 4(2):279--320, 1998.

\bibitem{MR1158626}
S.~Katz and D.~R. Morrison.
\newblock Gorenstein threefold singularities with small resolutions via
  invariant theory for {W}eyl groups.
\newblock {\em J. Algebraic Geom.}, 1(3):449--530, 1992.

\bibitem{MR992334}
P.~B. Kronheimer.
\newblock The construction of {ALE} spaces as hyper-{K}\"ahler quotients.
\newblock {\em J. Differential Geom.}, 29(3):665--683, 1989.

\bibitem{MR992335}
P.~B. Kronheimer.
\newblock A {T}orelli-type theorem for gravitational instantons.
\newblock {\em J. Differential Geom.}, 29(3):685--697, 1989.

\bibitem{MR710273}
U.~Lindstr{\"o}m and M.~Ro{\v{c}}ek.
\newblock Scalar tensor duality and {$N=1,\,2$} nonlinear {$\sigma $}-models.
\newblock {\em Nuclear Phys. B}, 222(2):285--308, 1983.

\bibitem{MR929144}
U.~Lindstr{\"o}m and M.~Ro{\v{c}}ek.
\newblock New hyper-{K}\"ahler metrics and new supermultiplets.
\newblock {\em Comm. Math. Phys.}, 115(1):21--29, 1988.

\bibitem{Lindstrom:2008gs}
U.~Lindstr{\"o}m and M.~Ro{\v{c}}ek.
\newblock {Properties of hyperk{\"a}hler manifolds and their twistor spaces}.
\newblock {\em Commun. Math. Phys.}, 293:257--278, 2010.

\bibitem{Lindstrom:1999pz}
U.~Lindstrom, M.~Ro\v{c}ek, and R.~von Unge.
\newblock {HyperK\"ahler quotients and algebraic curves}.
\newblock {\em JHEP}, 01:022, 2000.

\bibitem{MR1873248}
T.-T. Lu and S.-H. Shiou.
\newblock Inverses of {$2\times 2$} block matrices.
\newblock {\em Comput. Math. Appl.}, 43(1-2):119--129, 2002.

\bibitem{Majorana1932}
E.~Majorana.
\newblock Atomi orientati in campo magnetico variabile.
\newblock {\em Il Nuovo Cimento}, 9(2):43--50, 1932.

\bibitem{MR2778451}
V.~Minerbe.
\newblock On the asymptotic geometry of gravitational instantons.
\newblock {\em Ann. Sci. \'Ec. Norm. Sup\'er. (4)}, 43(6):883--924, 2010.

\bibitem{MR953820}
H.~Pedersen and Y.~S. Poon.
\newblock Hyper-{K}\"ahler metrics and a generalization of the {B}ogomolny
  equations.
\newblock {\em Comm. Math. Phys.}, 117(4):569--580, 1988.

\bibitem{MR1048125}
R.~Penrose.
\newblock {\em The emperor's new mind}.
\newblock The Clarendon Press, Oxford University Press, New York, 1989.

\bibitem{MR1865778}
R.~Penrose.
\newblock {\em Shadows of the mind}.
\newblock Oxford University Press, Oxford, 1994.

\bibitem{MR664330}
S.~Salamon.
\newblock Quaternionic {K}\"ahler manifolds.
\newblock {\em Invent. Math.}, 67(1):143--171, 1982.

\bibitem{Seiberg:1996nz}
N.~Seiberg and E.~Witten.
\newblock {Gauge dynamics and compactification to three dimensions}.
\newblock In {\em {The mathematical beauty of physics: A memorial volume for
  Claude Itzykson. Proceedings, Conference, Saclay, France, June 5-7, 1996}},
  pages 333--366, 1996.

\bibitem{Sen:1997kz}
A.~Sen.
\newblock {A Note on enhanced gauge symmetries in M and string theory}.
\newblock {\em JHEP}, 09:001, 1997.

\end{thebibliography}
\bibliographystyle{abbrv}

\end{document}